\documentclass{amsart}
\usepackage{amssymb}

\numberwithin{equation}{section}
\input{tcilatex}

\begin{document}
\title[Minimizers for the p-area]{Variations of generalized area functionals
and p-area minimizers of bounded variation in the Heisenberg group}
\author{Jih-Hsin Cheng}
\address[Cheng]{ Institute of Mathematics, Academia Sinica, Taipei and
National Center for Theoretical Sciences, Taipei Office, Taiwan, R.O.C.}
\email{cheng@math.sinica.edu.tw}
\urladdr{}
\thanks{}
\author{Jenn-Fang Hwang}
\address[Hwang]{ Institute of Mathematics, Academia Sinica, Taipei, Taiwan,
R.O.C.}
\email{majfh@math.sinica.edu.tw}
\urladdr{}
\thanks{}
\subjclass{Primary: 35L80; Secondary: 35J70, 32V20, 53A10, 49Q10.}
\keywords{Minimizer, p-area, BV, Heisenberg group, first variation, second
variation.}
\thanks{}

\begin{abstract}
We prove the existence of a continuous $BV$ minimizer with $C^{0}$ boundary
value for the $p$-area (pseudohermitian or horizontal area) in a
parabolically convex bounded domain. We extend the domain of the area
functional from $BV$ functions to vector-valued measures. Our main purpose
is to study the first and second variations of such a generalized area
functional including the contribution of the singular part. By giving
examples in Riemannian and pseudohermitian geometries, we illustrate several
known results in a unified way. We show the contribution of the singular
curve in the first and second variations of the $p$-area for a surface in an
arbitrary pseudohermitian $3$-manifold.
\end{abstract}

\maketitle

\bigskip

%$^{1}$Institute of Mathematics, Academia Sinica, Nankang, Taipei 11529,
%Taiwan, R.O.C. \and Jenn-Fang Hwang$^{2}$ \\
%EndAName
%$^{2}$Institute of Mathematics, Academia Sinica, Nankang, Taipei 11529,
%Taiwan, R.O.C.}
%\subjclass{}
%\email{cheng@math.sinica.edu.tw\\
%$\ \ \ \ \ ^{\ast }$Research supported in part by the NSC of Taiwan}
%\email{majfh@math.sinica.edu.tw}
%\subjclass{Primary 35L80; Secondary ??}
%\keywords{ P-area, p-minimizer, bounded variation, Heisenberg group. }

%\begin{abstract}
%\end{abstract}

\section{\textbf{Introduction and statement of the results}}

In \cite{chy}, Paul Yang and the authors proved the existence of a Lipschitz
continuous ($p$-)minimizer with $C^{2,\alpha }$ boundary value for the $p$%
-area (or horizontal area) in the space $W^{1,1}$ and the uniqueness of $p$%
-minimizers in the space $W^{1,2}$ among other things. In this paper, we
will prove the existence of a continuous $BV$ minimizer with $C^{0}$
boundary value for the $p$-area in a parabolically convex bounded domain.
Recall that the $p$-area is a special case of a more general area functional:%
\begin{equation}
\mathcal{F}_{H}(u):=\int_{\Omega }(|\nabla u+\vec{F}|+Hu)d^{m}x.
\label{eqn1.0}
\end{equation}

\noindent where $\Omega \subset R^{m}$ is a bounded domain, $u$ $\in $ $%
W^{1,1}(\Omega ),$ $\vec{F}$ is an $L^{1}$ vector field on $\Omega ,$ $H$ $%
\in $ $L^{\infty }(\Omega ),$ and $d^{m}x$ $:=$ $dx_{1}\wedge dx_{2}\wedge
...\wedge dx_{m}$ denotes the Euclidean volume form or the Lebesgue measure.
We often denote $\mathcal{F}_{H}$ by $\mathcal{F}$ for the case of $H$ $=$ $%
0:$%
\begin{equation}
\mathcal{F}(u):=\int_{\Omega }|\nabla u+\vec{F}|d^{m}x.  \label{eqn1.1}
\end{equation}%
\noindent $\mathcal{F(\cdot )}$ is called the $p$-area (of the graph defined
by $u$ over $\Omega )$ if $\vec{F}$ $=$ $-\vec{X}^{\ast }$ where $\vec{X}%
^{\ast }$ $=$ $(x_{1^{\prime }},$ $-x_{1},$ $x_{2^{\prime }},$ $-x_{2},...,$ 
$x_{n^{\prime }},$ $-x_{n}),$ $m$ $=$ $2n$ (see \cite{chmy}). In the case of
a graph $\Sigma $ over the $R^{2n}$-hyperplane in the Heisenberg group, the
above definition of $p$-area coincides with those given in \cite{CDG}, \cite%
{DGN}, and \cite{Pau}. In particular these notions, especially in the
framework of geometric measure theory, have been used to study existence or
regularity properties of minimizers for the relative perimeter or
extremizers of isoperimetric inequalities (see, e.g., \cite{DGN}, \cite{Ga}, 
\cite{LM}, \cite{LR}, \cite{MR}, \cite{Pan}, \cite{Pau2}, \cite{Ri}, \cite%
{Ri2}, \cite{CC}). The $p$-area can also be identified with the $2n+1$%
-dimensional spherical Hausdorff measure of $\Sigma $ (see, e.g., \cite{Ba}, 
\cite{FSS}). Some authors take the viewpoint of so called intrinsic graphs
and obtained interesting results (see, e.g., \cite{FSS}, \cite{ASCV}, and 
\cite{BSC} which relates distributional solutions of Burgers' equation to
intrinsic regular graphs). Starting from the work \cite{chmy} (see also \cite%
{ch}), the subject was studied from the viewpoint of partial differential
equations and that of differential geometry (see \cite{chy}, \cite{chy1}, 
\cite{chmy2}; the term $p$-minimal is used since this is the notion of
minimal surfaces in pseudohermitian geometry; "$p$" stands for
"pseudohermitian"). In \cite{chy}, one studied the situation for $u\in
W^{1,1}.$ To extend the domain of $\mathcal{F}$ to the space of $BV$
functions, we define the total variation of a function $u\in L^{1}(\Omega )$
by

\begin{equation}
\int_{\Omega }|Du+\vec{F}d^{m}x|\text{ :}=\sup \{\int_{\Omega }(-u\func{div}%
\vec{\phi}+\vec{F}\cdot \vec{\phi})d^{m}x\mid \vec{\phi}\in C_{0}^{1}(\Omega
),|\vec{\phi}|\leq 1\}.  \label{eqn1.2}
\end{equation}

\noindent Let $BV_{\vec{F}}(\Omega )$ denote the space of $u$ $\in $ $%
L^{1}(\Omega )$ such that the total variation $\int_{\Omega }|Du+\vec{F}%
d^{m}x|$ $<$ $\infty .$ In this case, the notation $Du$ (viewed as the
gradient of $u$ in the distributional sense) is in fact a vector-valued
Radon signed measure (see Remark 1.5 on page 5 in \cite{Gi}) and $|Du+\vec{F}%
d^{m}x|$ is the total variation measure of the measure $Du+\vec{F}d^{m}x$
(see the first paragraph of Section 3 for more details). When $u$ $\in $ $%
W^{1,1}(\Omega ),$ we use $\nabla u$ to denote the gradient of $u.$ Note
that $BV_{\vec{F}}(\Omega )$ is reduced to the usual space of $BV$
functions, denoted by $BV(\Omega ),$ for $\vec{F}=\vec{0}.$ Moreover, if $%
\vec{F}$ $\in $ $L^{1}(\Omega ),$ it is easy to see that $u\in BV_{\vec{F}%
}(\Omega )$ if and only if $u\in BV(\Omega ).$ For $u$ $\in $ $W^{1,1}$ (\ref%
{eqn1.2}) is the same as the one in the usual sense (in which we write $Du$ $%
=$ ($\nabla u)d^{m}x)$.

We need to require the following condition on $\vec{F}$ (say, $\in $ $%
C^{1}): $

\begin{equation}
\partial _{K}F_{I}=\partial _{I}f_{K},\ \ I,K=1,...,m  \label{eqn1.3}
\end{equation}

\noindent for $C^{1}$-smooth functions $f_{K}$'s. Denote the coordinates of $%
R^{m}$ by $x_{1},$ $x_{2},$ $...,$ $x_{m}.$ We call a coordinate system
orthonormal if it is obtained by a translation and a rotation from $x_{1},$ $%
x_{2},$ $...,$ $x_{m}.$ We recall (\cite{chy}) the definition of a certain
notion of convexity for $\Omega $ as follows.

\bigskip

\textbf{Definition 1.1.} We call $\Omega $ $\subset $ $R^{m}$ parabolically
convex if for any $p\in \partial \Omega ,$ there exists an orthonormal
coordinate system $(\tilde{x}_{1},$ $\tilde{x}_{2},$ $...,$ $\tilde{x}_{m})$
with the origin at $p$ and $\Omega \subset \{a\tilde{x}_{1}^{2}-\tilde{x}%
_{2}<0\}$ where $a>0$ is independent of $p.$

\bigskip

Note that a $C^{2}$-smooth bounded domain with the positively curved
(positive principal curvatures) boundary is parabolically convex. On the
other hand, a parabolically convex domain can be nonsmooth as shown by the
following example: a planar domain defined by 
\begin{equation*}
-\sqrt{3}<x<\sqrt{3},\text{ }-\sqrt{4-x^{2}}+1<y<\sqrt{4-x^{2}}-1.
\end{equation*}%
\textit{\noindent }For a vector field $\vec{G}$ $=$ $(g_{1},g_{2},...,g_{2n})
$ on $\Omega \subset R^{2n},$ we define $\vec{G}^{\ast }$ $:=$ $(g_{2},$ $%
-g_{1},$ $g_{4},$ $-g_{3},$ $...,$ $g_{2n},$ $-g_{2n-1}).$

\bigskip

\textbf{Theorem A. }\textit{Let }$\Omega $ \textit{be a parabolically} 
\textit{convex bounded domain in }$R^{2n}$ \textit{with }$\partial \Omega
\in C^{2,\alpha }$\textit{\ }$(0<\alpha <1)$\textit{. Suppose }$\vec{F}$%
\textit{\ }$\in $\textit{\ }$C^{1,\alpha }(\bar{\Omega})$\textit{\ satisfies
the condition (\ref{eqn1.3}) for }$C^{1,\alpha }$\textit{-smooth and bounded 
}$f_{K}$\textit{'s in }$\Omega $ \textit{and }$\func{div}\vec{F}^{\ast }$%
\textit{\ }$>$\textit{\ }$0.$ \textit{Let }$\varphi \in C^{0}(\partial
\Omega ).$\textit{\ Then there exists }$u\in C^{0}(\bar{\Omega})\cap
BV(\Omega )$\textit{\ such that }$u=\varphi $\textit{\ on }$\partial \Omega $%
\textit{\ and}

\begin{equation}
\int_{\Omega }|Du+\vec{F}d^{m}x|\text{ }\leq \int_{\Omega }|Dv+\vec{F}d^{m}x|%
\text{ }  \label{eqn1.4}
\end{equation}%
\textit{\noindent for all }$v\in C^{0}(\bar{\Omega})\cap BV(\Omega )$\textit{%
\ with }$v=\varphi $\textit{\ on }$\partial \Omega .$

\bigskip

We remark that $\vec{F}$ $=$ $-\vec{X}^{\ast }$ satisfies the assumption in
Theorem A. The idea of the proof for Theorem A goes as follows. We
approximate $\varphi $ by $C^{2,\alpha }$-smooth functions and apply Theorem
A in \cite{chy} to get approximating Lipschitz continuous minimizers. These
minimizers will converge uniformly to a continuous function $u$ by the
comparison principle (Theorem C in \cite{chy} ). Then we show that $u$ is a $%
BV$ function and a minimizer in $C^{0}(\bar{\Omega})\cap BV(\Omega )$ by
some extra work$.$

On the other hand, F. Serra Cassano and D. Vittone in a recent paper (\cite%
{SCV}) study this problem for more general domains. Let $\Omega $ $\subset $ 
$R^{2n}$ be a bounded domain with Lipschitz regular boundary. They show the
functional%
\begin{equation}
u\in BV(\Omega )\rightarrow \int_{\Omega }|Du-\vec{X}^{\ast
}d^{2n}x|+\int_{\partial \Omega }\mid u|_{\partial \Omega }-\varphi \mid
d\sigma  \label{eqn1.4.0}
\end{equation}
\textit{\noindent }attains its minimum, where,\textit{\ }for $u$ $\in $ $%
BV(\Omega ),$ the trace $u|_{\partial \Omega }$ exists and lies in $%
L^{1}(\partial \Omega )$ by Theorem 2.10 in \cite{Gi}, $\varphi $ $\in $ $%
L^{1}(\partial \Omega )$ is given, and $d\sigma $ denotes the standard
boundary measure. Moreover, there holds

\begin{eqnarray*}
&&\inf \{\int_{\Omega }|Du-\vec{X}^{\ast }d^{2n}x| \text{:} u\in BV(\Omega
),u|_{\partial \Omega }=\varphi \} \\
&=& \min \{\int_{\Omega }|Du-\vec{X}^{\ast }d^{2n}x|+\int_{\partial \Omega
}\mid u|_{\partial \Omega }-\varphi \mid d\sigma :u\in BV(\Omega )\}
\end{eqnarray*}

\noindent (see Theorem 1.4 in \cite{SCV}).

Although the $BV$ minimizers $\tilde{u}$ for (\ref{eqn1.4.0}) exist, the
trace $\tilde{u}|_{\partial \Omega }$ may not equal $\varphi .$ The $BV$
minimizers for $\int_{\Omega }|Du-\vec{X}^{\ast }d^{2n}x|$ with given (even
smooth) boundary value $\varphi $ may not exist in general either for
nonconvex domains as shown in Example 3.6 of \cite{SCV}. In fact, consider $%
\Omega $ $:=$ $\{1$ $<$ $\sqrt{x^{2}+y^{2}}$ $<$ $2\}$ $\subset $ $R^{2}.$
Take the boundary value $\varphi $ $=$ $0$ on $\sqrt{x^{2}+y^{2}}$ $=$ $2$
while $\varphi $ $=$ $C$ on $\sqrt{x^{2}+y^{2}}$ $=$ $1.$ Then there admits
no minimizer for $\int_{\Omega }|Du-\vec{X}^{\ast }d^{2n}x|$ with $%
u|_{\partial \Omega }$ $=$ $\varphi $ when $C$ is large enough (see \cite%
{SCV} for more details). The original idea comes from \cite{Fi} in which R.
Finn gave examples of nonexistence for the Dirichlet problem of (Euclidean)
minimal surface equation.

After we have $BV$ minimizers, we consider the variations of $\mathcal{F}$
on $BV$ functions. Since $\mathcal{F}$ is only convex, but not strongly
convex, this causes much trouble. Besides the trouble that $BV$ functions
cause, we still have trouble even for $C^{\infty }$-smooth functions. For
instance, let $\vec{F}$ $=$ $\vec{0},$ then $u$ $\equiv $ $0$ is the
minimizer for $\mathcal{F(}u\mathcal{)}$ $=$ $\int_{\Omega }|\nabla u|d^{m}x$%
. Compute the first variation at $u$ $\equiv $ $0$:%
\begin{eqnarray}
&&\lim_{\varepsilon \rightarrow 0}\frac{\mathcal{F(}0+\varepsilon \varphi )-%
\mathcal{F(}0)}{\varepsilon }  \label{eqn1.4.1} \\
&=&\lim_{\varepsilon \rightarrow 0}\frac{|\varepsilon |}{\varepsilon }%
\int_{\Omega }|\nabla \varphi |d^{m}x.  \notag
\end{eqnarray}%
\textit{\noindent }from which\textit{\ }we learn that only left limit or
right limit exists. However, we can still deal with the second variation of $%
\mathcal{F}$ (see Theorem C). Previously in \cite{chmy} the second variation
of $\mathcal{F}$ was studied only for $C^{2}$ smooth functions and away from
the singular set $S_{\vec{F}}(u)$ ($:=$ $\{p$ $\in $ $\Omega $ $|$ $\nabla u+%
\vec{F}$ $=$ $0$ at $p\}).$ But whether $H_{m}(S_{\vec{F}}(u))$ ($m$ $=$ $%
\dim \Omega ),$ the $m$-th dimensional Hausdorff measure of $S_{\vec{F}}(u),$
vanishes is a problem. In the case of least gradient ($\vec{F}$ $=$ $0),$ $%
H_{m}(S_{\vec{F}}(u))$ may not be zero.

In the case of $p$-area, $\vec{F}$ $=$ $-\vec{X}^{\ast }$ where $\vec{X}%
^{\ast }$ $=$ $(x_{1^{\prime }},$ $-x_{1},$ $x_{2^{\prime }},$ $-x_{2},...,$ 
$x_{n^{\prime }},$ $-x_{n}),$ $m$ $=$ $2n,$ for $u$ $\in $ $C^{2}(\Omega ),$
we have $H_{m}(S_{\vec{F}}(u))$ $=$ $0.$ But for $u$ $\in $ $W^{1,1}(\Omega
),$ $H_{m}(S_{\vec{F}}(u))$ may be larger than zero (see \cite{Ba}). For $u$ 
$\in $ $BV(\Omega ),$ we write 
\begin{equation*}
Du=(\nabla u)d^{m}x+dv_{s}^{d^{m}x},\text{ }d^{m}x\perp dv_{s}^{d^{m}x}
\end{equation*}%
\textit{\noindent }where\textit{\ }$\nabla u$ $\in $ $L^{1}$ with respect to 
$d^{m}x$. Suppose $d^{m}x$ ($dv_{s}^{d^{m}x},$ resp.) is concentrated on $%
\Omega _{1}$ ($\Omega _{2},$ resp.) where $\Omega _{1}$ $\cap $ $\Omega _{2}$
$=$ $\varnothing ,$ $\Omega _{1}$ $\cup $ $\Omega _{2}$ $=$ $\Omega .$ Note
that $H_{m}(\Omega _{2})$ $=$ $0.$ Define $S_{\vec{F}}(u)$ $:=$ $\{p$ $\in $ 
$\Omega _{1}$ $|$ $\nabla u+\vec{F}$ $=$ $0$ at $p\}.$ Now whether $H_{m}(S_{%
\vec{F}}(u))$ $=$ $0$ ($m$ even) for a $BV$ minimizer $u$ for the $p$-area
in general is still an open problem. So we cannot neglect the role of $S_{%
\vec{F}}(u).$ One of the purposes of this paper is to study the second
variation of $\mathcal{F}$ not avoiding $S_{\vec{F}}(u)$ even if $H_{m}(S_{%
\vec{F}}(u))$ $\neq $ $0.$

The idea of computing the first and second variations is to extend the
domain of $\mathcal{F(\cdot )}$ from $BV$ functions to vector-valued
measures. Then making use of the Radon-Nikodym theorem, we can easily obtain
the formulas of first and second variations, which include the effect of the
singular set.

Let $E$ be a $C^{\infty }$-smooth Riemannian vector bundle over a $C^{\infty
}$-smooth manifold $X$. Let $d\mu ,$ $d\nu $ be two $E$-valued (Radon
signed) measures on $X$. Let $d\mu _{\varepsilon }$ $=$ $d\mu $ $+$ $%
\varepsilon d\nu $ for $\varepsilon $ $\in $ $R.$ Define $\mathcal{F(}d\mu
_{\varepsilon })$ by%
\begin{equation*}
\mathcal{F(}d\mu _{\varepsilon }):=\int_{X}|d\mu _{\varepsilon }|.
\end{equation*}%
\textit{\noindent }(see (\ref{eqn3.1}) with $\Omega $ replaced by $X)$
Denote $\mathcal{F(}d\mu _{\varepsilon })$ by $\mathcal{F(\varepsilon )}$
for simplicity. Throughout this paper we assume that both $d\mu $ and $d\nu $
are bounded in the sense that $|d\mu |$ and $|d\nu |$ are integrable over $%
X. $ By the (extended) Radon-Nikodym theorem we can write

\begin{eqnarray*}
d\mu _{\varepsilon } &=&N_{\varepsilon }|d\mu _{\varepsilon }|, \\
d\nu  &=&A_{\varepsilon }|d\mu _{\varepsilon }|+d\nu _{s}^{\varepsilon },%
\text{ }d\nu _{s}^{\varepsilon }\perp |d\mu _{\varepsilon }|
\end{eqnarray*}%
\textit{\noindent }where $N_{\varepsilon }$, $A_{\varepsilon }$ $\in $ $%
L^{1}(|d\mu _{\varepsilon }|)$ with $|N_{\varepsilon }|$ $=$ $1$ (cf. (\ref%
{eqn3.2}))$.$ Recall that $\mathcal{F}^{\prime }\mathcal{(\varepsilon }_{1}%
\mathcal{\pm )}$ $:=$ $\lim_{\varepsilon _{2}\rightarrow \varepsilon _{1}\pm
}\frac{\mathcal{F(\varepsilon }_{2}\mathcal{)-F(\varepsilon }_{1}\mathcal{)}%
}{\varepsilon _{2}-\varepsilon _{1}}.$ We have the following first variation
formula.

\bigskip

\textbf{Theorem B.} \textit{Suppose }$d\mu $\textit{\ and }$d\nu $\textit{\
are bounded. Then }$\mathcal{F(\varepsilon )}$ \textit{is Lipschitz
continuous in }$\varepsilon $\textit{\ and} \textit{there holds}

\begin{equation}
\mathcal{F}^{\prime }\mathcal{(\varepsilon }_{1}\mathcal{\pm )=}%
\int_{X}N_{\varepsilon _{1}}\cdot A_{\varepsilon _{1}}|d\mu _{\varepsilon
_{1}}|\pm \int_{X}|d\nu _{s}^{\varepsilon _{1}}|.  \label{eqn1.5}
\end{equation}

\bigskip

Let $u$ $\in $ $BV(\Omega )$ where $\Omega $ $\subset $ $R^{m}$ is a bounded
domain with Lipschitz regular boundary. Define%
\begin{equation*}
\mathcal{\tilde{F}}_{H}(u):=\int_{\Omega }|Du+\vec{F}d^{m}x|+\int_{\Omega }Hu%
\text{ }d^{m}x
\end{equation*}%
\textit{\noindent }where the first term on the right side of the equality
makes sense by (\ref{eqn1.2}). Recall that $\vec{F}$ is an $L^{1}$ vector
field on $\Omega $ and $H$ $\in $ $L^{\infty }(\Omega ).$ Recall that for $u$
$\in $ $BV(\Omega ),$ the trace $u|_{\partial \Omega }$ exists and lies in $%
L^{1}(\partial \Omega )$ by Theorem 2.10 in \cite{Gi}.

\bigskip

\textbf{Definition 1.2.} Suppose $u$ $\in $ $BV(\Omega )$ with $u|_{\partial
\Omega }$ $=$ $\psi $. If for all $\varphi $ $\in $ $BV(\Omega )$ with $%
\varphi |_{\partial \Omega }$ $=$ $0,$ there holds 
\begin{equation*}
\mathcal{\tilde{F}}_{H}(u)\leq \mathcal{\tilde{F}}_{H}(u+\varphi )
\end{equation*}%
\textit{\noindent }Then we call $u$ a minimizer for $\mathcal{\tilde{F}}_{H}$
with the boundary value (trace) $\psi .$

\bigskip

Denote $\mathcal{\tilde{F}}_{H}(u+\varepsilon \varphi )$ by $\mathcal{\tilde{%
F}}_{H}(\varepsilon )$.\textit{\ }We can then have the following necessary
conditions for $u$ $\in $ $BV(\Omega )$ to be a minimizer.

\bigskip

\textbf{Corollary B}$^{\prime }$\textbf{.} \textit{Let }$\Omega $\textit{\ }$%
\subset $\textit{\ }$R^{m}$\textit{\ be a bounded domain with Lipschitz
regular boundary.} \textit{Suppose }$u$ $\in $ $BV(\Omega )$ \textit{is a
minimizer for }$\mathcal{\tilde{F}}_{H}$ \textit{with }$u|_{\partial \Omega
} $ \textit{= }$\psi $ $\in $ $L^{1}(\partial \Omega )$\textit{. Then} 
\textit{there hold}%
\begin{equation}
\mathcal{\tilde{F}}_{H}^{\prime }\mathcal{(}0+\mathcal{)=}\int_{\Omega
}N_{0}\cdot A_{0}|d\mu |+\int_{\Omega }|d\nu _{s}^{0}|+\int_{\Omega
}H\varphi \text{ }d^{m}x\geq 0  \label{eqn1.5.0}
\end{equation}%
\textit{\noindent and}%
\begin{equation}
\mathcal{\tilde{F}}_{H}^{\prime }\mathcal{(}0-\mathcal{)=}\int_{\Omega
}N_{0}\cdot A_{0}|d\mu |-\int_{\Omega }|d\nu _{s}^{0}|+\int_{\Omega
}H\varphi \text{ }d^{m}x\leq 0.  \label{eqn1.5.1}
\end{equation}

\bigskip

We remark that Corollary B$^{\prime }$ generalizes Theorem 3.3 in \cite{chy}%
, where $N(u)$ $=$ $N_{0},$ ($\nabla \varphi )d^{m}x$ $=$ $A_{0}|d\mu |$ on $%
\Omega \backslash S_{\vec{F}}(u),$ and \TEXTsymbol{\vert}$\nabla \varphi
|d^{m}x$ $=$ $|d\nu _{s}^{0}|$ on $S_{\vec{F}}(u).$ Here $A_{0}$ $=$ $\frac{%
\nabla \varphi }{|\nabla u+\vec{F}|},$ $|d\mu |$ $=$ $|\nabla u+\vec{F}%
|d^{m}x,$ and $S_{\vec{F}}(u)$ :$=$ $\{\nabla u+\vec{F}$ $=$ $0\}$ (cf.
Example 3.2)$.$ Also note that (\ref{eqn1.5.1}) corresponds to (3.12) in 
\cite{chy} with $\varphi $ replaced by $-\varphi .$

The singular term $\pm \int_{X}|d\nu _{s}^{0}|$ in (\ref{eqn1.5.0}) and (\ref%
{eqn1.5.1}) is not removable in general. The simplest example is that at the
minimizer $u$ $\equiv $ $0$ for the least gradient energy functional $\int
|\nabla u|d^{m}x,$ we have $\mathcal{F}^{\prime }(0\pm )$ $=$ $\pm \int
|d\nu _{s}^{0}|$ $=$ $\pm \int $\TEXTsymbol{\vert}$\nabla \varphi |d^{m}x$
over $S_{\vec{F}}(u)$ $=$ $\Omega $ as shown in (\ref{eqn1.4.1}).

There are at most countably many $\varepsilon $'s such that $|d\nu
_{s}^{\varepsilon }|$ $\neq $ $0.$ We call $\varepsilon $ regular if $|d\nu
_{s}^{\varepsilon }|$ $=$ $0.$ For regular $\varepsilon $ we have $\mathcal{F%
}^{\prime }\mathcal{(\varepsilon }+\mathcal{)}$ $=$ $\mathcal{F}^{\prime }%
\mathcal{(\varepsilon }-\mathcal{)}$ $=$ $\mathcal{F}^{\prime }\mathcal{%
(\varepsilon )}$ and $\mathcal{F}^{\prime }\mathcal{(\varepsilon )}$ is an
increasing function of $\varepsilon $ (see Proposition 3.3)$.$ Write $%
\mathcal{F}_{+}^{\prime }\mathcal{(\varepsilon )}$ ($\mathcal{F}_{-}^{\prime
}\mathcal{(\varepsilon )}$, respectively) for $\mathcal{F}^{\prime }\mathcal{%
(\varepsilon }+\mathcal{)}$ ($\mathcal{F}^{\prime }\mathcal{(\varepsilon }-%
\mathcal{)}$, respectively). In Section 3 we also study the left and right
continuity of $\mathcal{F}_{+}^{\prime }$ and $\mathcal{F}_{-}^{\prime }$
(see Proposition 3.4)$.$ We give area functionals in Riemannian and
pseudohermitian geometries as examples to illustrate (\ref{eqn1.5}). For a $%
p $-area stationary surface in an arbitrary pseudohermitian 3-manifold, we
obtain the "incident angle = reflected angle" condition on the singular
curve (see (\ref{eqn3.30})). The result extends previous ones in the
Heisenberg group (\cite{chy}, \cite{Ri2}).

In Section 4 we discuss the second derivative of $\mathcal{F(\varepsilon )}$%
. We compute the first derivatives of $\mathcal{F}_{+}^{\prime }$ and $%
\mathcal{F}_{-}^{\prime }$ in various situations.

\bigskip

\textbf{Theorem C.} \textit{Suppose }$d\mu $\textit{\ and }$d\nu $\textit{\
are bounded, and }$|A_{\varepsilon _{1}}|^{2}$ $\in $ $L^{1}(X,|d\mu
_{\varepsilon _{1}}|).$ \textit{Then} \textit{(1) For\ }$\varepsilon _{1}$%
\textit{\ regular, there holds}

\begin{equation}
\lim_{\varepsilon _{2}\rightarrow \varepsilon _{1},\text{ }\varepsilon _{2}%
\text{ regular}}\frac{\mathcal{F}^{\prime }\mathcal{(}\varepsilon _{2})-%
\mathcal{F}^{\prime }\mathcal{(}\varepsilon _{1})}{\varepsilon
_{2}-\varepsilon _{1}}=\int_{X}\{|A_{\varepsilon _{1}}|^{2}-|(A_{\varepsilon
_{1}}\cdot N_{\varepsilon _{1}})|^{2}\}|d\mu _{\varepsilon _{1}}|\text{ \ }%
(\geq 0).  \label{eqn1.6}
\end{equation}%
\textit{\noindent (2) For }$\varepsilon _{1}$\textit{\ arbitrary, there holds%
}

\begin{eqnarray}
\lim_{\varepsilon _{2}\rightarrow \varepsilon _{1}+}\frac{\mathcal{F}_{\pm
}^{\prime }\mathcal{(}\varepsilon _{2})-\mathcal{F}_{+}^{\prime }\mathcal{(}%
\varepsilon _{1})}{\varepsilon _{2}-\varepsilon _{1}} &=&\lim_{\varepsilon
_{2}\rightarrow \varepsilon _{1}-}\frac{\mathcal{F}_{\pm }^{\prime }\mathcal{%
(}\varepsilon _{2})-\mathcal{F}_{-}^{\prime }\mathcal{(}\varepsilon _{1})}{%
\varepsilon _{2}-\varepsilon _{1}}  \label{eqn1.7} \\
&=&\int_{X}\{|A_{\varepsilon _{1}}|^{2}-|(A_{\varepsilon _{1}}\cdot
N_{\varepsilon _{1}})|^{2}\}|d\mu _{\varepsilon _{1}}|\text{ \ }(\geq 0). 
\notag
\end{eqnarray}

\bigskip

Observe that $\mathcal{F}_{-}^{\prime }(\varepsilon _{1})$ may be strictly
less than $\mathcal{F}_{+}^{\prime }\mathcal{(}\varepsilon _{1})$ (roughly
speaking, $\mathcal{F}^{\prime }$ is not continuous and may have a jump at $%
\varepsilon _{1})$. Still we have not only the existence of the left
derivative of $\mathcal{F}_{-}^{\prime }$ and the right derivative of $%
\mathcal{F}_{+}^{\prime }$, but also the same value, i.e., $(\mathcal{F}%
_{-}^{\prime })_{-}^{\prime }(\varepsilon _{1})$ $=$ ($\mathcal{F}%
_{+}^{\prime })_{+}^{\prime }\mathcal{(}\varepsilon _{1})$ by (\ref{eqn1.7}%
). This is a very special property. Note that a convex function does not
have such a property in general. For instance, $f(x)$ $=$ $0$ for $x$ $\leq $
$0,$ $f(x)$ $=$ $x^{2}$ $+$ $x$ for $x$ $>$ $0.$ We can easily check that $%
f^{\prime }$ has a jump at $x$ $=$ $0.$ On the other hand, we compute $%
f^{\prime \prime }(x)$ = $0$ for $x$ $<$ $0$ while $f^{\prime \prime }(x)$ = 
$2$ for $x$ $>$ $0.$

A fundamental formula in deducing the second variation of $\mathcal{F}$ is (%
\ref{eqn3.11}) (for \TEXTsymbol{\vert}$d\mu _{\varepsilon _{1}}|$ $\ll $ 
\TEXTsymbol{\vert}$d\mu _{\varepsilon _{2}}|$ $\ll $ \TEXTsymbol{\vert}$d\mu
_{\varepsilon _{1}}|$) in Section 3:%
\begin{equation*}
(N_{\varepsilon _{2}}-N_{\varepsilon _{1}})\cdot (d\mu _{\varepsilon
_{2}}-d\mu _{\varepsilon _{1}})=\frac{1}{2}|N_{\varepsilon
_{2}}-N_{\varepsilon _{1}}|^{2}(|d\mu _{\varepsilon _{2}}|+|d\mu
_{\varepsilon _{1}}|).
\end{equation*}%
\textit{\noindent }This formula generalizes (5.1) in \cite{chmy}:%
\begin{equation}
(N(u)-N(v))\cdot (\nabla u-\nabla v)=\frac{1}{2}|N(u)-N(v)|^{2}(|\nabla u-%
\vec{X}^{\ast }|+|\nabla v-\vec{X}^{\ast }|)  \label{1.7.1}
\end{equation}%
\textit{\noindent }for $u,$ $v$ $\in $ $C^{1}.$ The extension (\ref{eqn3.11}%
) includes the case of $BV$ functions. Also it holds for various geometries
including those of Euclidean and pseudohermitian minimal surfaces. See the
examples in Section 3 and the Appendix. Corresponding to (\ref{1.7.1}), for
the Riemannian mean curvature equation $\func{div}Tu$ $=$ $H$ in $R^{n},$
where $Tu$ :$=$ $\frac{\nabla u}{\sqrt{1+|\nabla u|^{2}}},$ we have the
following structural inequality:%
\begin{eqnarray*}
(Tu-Tv)\cdot (\nabla u-\nabla v) &\geq &\frac{1}{2}|Tu-Tv|^{2}(\sqrt{%
1+|\nabla u|^{2}}+\sqrt{1+|\nabla v|^{2}}) \\
&\geq &|Tu-Tv|^{2}.
\end{eqnarray*}%
\textit{\noindent }The\textit{\ }above inequality was discovered by
Miklyukov \cite{Mik}, Hwang \cite{Hw1}, and Collin-Krust \cite{CK}
independently. The proof in \cite{Hw1} was obtained through the help of
Shuh-Jye Chern who simplified the original proof of Hwang.

In Section 4 we give a proof of Theorem C and examples to illustrate (\ref%
{eqn1.6}). In particular we show that a $C^{2}$ area-stationary graph in a
flat ambient space in either Riemannian or pseudohermitian geometry has the
local area- minimizing property. This fact was proved individually for
different situations. For the 3-dimensional Heisenberg group, it was shown
by a calibration argument in \cite{chmy} for the nonsingular case. Later
Ritor\'{e} and Rosales (\cite{Ri2}) extended the result to the situation
having singularities. On the other hand, using (\ref{eqn1.6}) gives a
unified proof (see Example 4.1 and Example 4.2). Note that in (\cite{Ri2}),
we are in $C^{2}$-smooth category. The singular set has no contribution to
the second variation since its Lebesgue measure (in $R^{2n})$ vanishes
according to a result of Balogh (\cite{Ba}). Here Theorem C generalizes to
include the singular set contribution. On the other hand, we obtained
Balogh's result (for a $C^{2}$-smooth function) as Lemma 5.4 in \cite{chmy}
by a different argument (we used only elementary linear algebra and the
implicit function theorem in the proof). Later we generalized this result to
the situation of general $\vec{F}$ (see Theorem D in \cite{chy}).

When the ambient space is not flat, we know that the curvature appears in
the second variation formula and the second variation is no longer
nonnegative in general. This means that the way we vary by considering $%
|d\mu +\varepsilon dv|$ is not generic for nonflat ambient spaces. For a
variational vector field with support containing a singular curve, we
compute the second variation of the $p$-area for a stationary surface in
such a direction, and cook out the contribution of the singular curve (see (%
\ref{eqn4.25}); the computation was completed by Hung-Lin Chiu). Note that
in \cite{chmy} we have done such a computation for a variational vector
field with support away from the singular set.

In the Appendix, we define the notions of gradient and hypersurface area in
a general formulation unifying Riemannian and pseudohermitian (horizontal or
Heisenberg) structures for further development. In fact, these different
geometric structures on a differentiable manifold $M$ are better described
in a unified way by assigning a nonnegative inner product $<\cdot ,\cdot >$
on its cotangent bundle $T^{\ast }M$. The gradient $\nabla \varphi $ of a
smooth function $\varphi $ on $M$ with respect to these different geometric
structures can be expressed in a unified way as $\nabla \varphi $ :$=$ $%
G(d\varphi )$ where $G$ $:$ $T^{\ast }M$ $\rightarrow $ $TM$ is a natural
bundle morphism defined by $<G(\omega ),\eta >$ $=$ $<\omega ,\eta >$ for $%
\omega ,$ $\eta $ $\in $ $T^{\ast }M$ (cf. (\ref{H3}) and note that the
first $<\cdot ,\cdot >$ denotes the pairing between $TM$ and $T^{\ast }M).$
For instance, this $\nabla \varphi $ is nothing but the subgradient $\nabla
_{b}\varphi $ in the pseudohermitian case. The geometric information is
hidden in $G.$ If $\varphi $ is a defining function of a hypersurface $%
\Sigma $ $\subset $ $M,$ we can give a unified definition of area (or
volume) element $dv_{\Sigma }$ of $\Sigma $ as follows (cf. (\ref{H4})):%
\begin{equation*}
dv_{\Sigma }=\frac{d\varphi }{|d\varphi |}\text{ }\rfloor \text{ }dv_{M}
\end{equation*}%
\noindent where $dv_{M}$ is a volume form. This formula encodes the
Euclidean area element, the $p$- (or $H$-) area element for a graph or an
intrinsic graph in the Heisenberg group (see Example A.1, Example A.2, and
Example A.3, resp.; see also Examples A.4 and A.5 for a surface in a general
pseudohermitian 3-manifold). In particular, we recover the definition of
Ritor\'{e} and Rosales for the $p$- (or $H$-) area element (\cite{Ri2}). See
(\ref{H16}) in Example A.5 for more details.

We also derive a general formula for the mean curvature and give a number of
examples to illustrate it. See (\ref{H22}) and Examples A.6 and A.7. When
this work was being done, we received an interesting preprint (see \cite{SCV}%
) from Francesco Serra Cassano. In \cite{SCV}, the authors also studied the
existence (and local boundedness) of $BV$ minimizers for the $p$-area (of
what the authors call $t$-graphs and $X_{1}$-graphs). The definition (1.3)
that we use here is $S_{t}(u)$ on page 16 of \cite{SCV}. Also the boundary
value in \cite{SCV} is more general (see previous comments after Theorem A
for more details).

\bigskip

\textbf{Added in proof:} The authors were informed of papers \cite{HRR}, 
\cite{MSCV} in which the second variation of the $p$-area was also studied
and discussed.

\bigskip

\textbf{Acknowledgments.} The first author's research was supported in part
by NSC 97-2115-M-001-016-MY3. He would also like to thank the National
Center for Theoretical Sciences, Taipei Office\ for sponsoring the Workshop
on Mean Curvature Equation in Heisenberg Geometry held December 18-19, 2010
at the Academia Sinica, Taipei. The second author's research was supported
in part by NSC 97-2115-M-001-005-MY3.

\bigskip

\section{\textbf{Existence and} \textbf{proof of Theorem A}}

Take $\varphi _{j}^{-},$ $\varphi _{k}^{+}$ $\in $ $C^{\infty }(\bar{\Omega}%
) $ such that $\varphi _{j}^{-}$ ($\varphi _{k}^{+},$ respectively)
increasingly (decreasingly, respectively) approaches $\varphi $ in $C^{0}$%
-norm on $\partial \Omega .$ By Theorem A in \cite{chy}, we can find
Lipschitz continuous minimizers (for $\mathcal{F}(\cdot ))$ $u_{j}^{-},$ $%
u_{k}^{+}$ such that $u_{j}^{-}$ $=$ $\varphi _{j}^{-}$ and $u_{k}^{+}$ $=$ $%
\varphi _{k}^{+}$ on $\partial \Omega .$ It follows from the maximum
principle (Theorem C in \cite{chy}; here the condition $\func{div}\vec{F}%
^{\ast }$ $>$ $0$ is used) that

\begin{eqnarray}
0 &\leq &u_{j_{1}}^{-}-u_{j_{2}}^{-}\leq \Vert \varphi _{j_{1}}^{-}-\varphi
_{j_{2}}^{-}\Vert _{C^{0}(\partial \Omega )},  \label{eqn2.1} \\
0 &\leq &u_{k_{1}}^{+}-u_{k_{2}}^{+}\leq \Vert \varphi _{k_{1}}^{+}-\varphi
_{k_{2}}^{+}\Vert _{C^{0}(\partial \Omega )},\text{ and}  \notag \\
0 &\leq &u_{k}^{+}-u_{j}^{-}\leq \Vert \varphi _{k}^{+}-\varphi
_{j}^{-}\Vert _{C^{0}(\partial \Omega )}  \notag
\end{eqnarray}

\noindent in $\Omega $ for $j_{1}\geq j_{2,}$ $k_{1}\leq k_{2}$ (note that
if $u$ is a solution or a minimizer, so is $u$ $+$ $a$ $constant$).Therefore
in view of (\ref{eqn2.1}) $u_{j}^{-}$ increasingly and $u_{k}^{+}$
decreasingly converge to the same limit $u$ $\in $ $C^{0}(\bar{\Omega})$
such that $u$ $=$ $\varphi $ on $\partial \Omega .$

\bigskip

\textbf{Lemma 2.1}. \textit{Let }$w_{j}$\textit{\ }$\in $\textit{\ }$%
L^{1}(\Omega )$\textit{\ }$\cap $\textit{\ }$BV_{\vec{F}}(\Omega ),$\textit{%
\ }$w$\textit{\ }$\in $\textit{\ }$L^{1}(\Omega ).$\textit{\ Suppose }$w_{j}$%
\textit{\ }$\rightarrow $\textit{\ }$w$\textit{\ in }$L^{1}.$\textit{\ Then}%
\begin{equation}
\int_{\Omega }|Dw+\vec{F}d^{m}x|\leq \lim \inf_{j\rightarrow \infty
}\int_{\Omega }|Dw_{j}+\vec{F}d^{m}x|.  \label{eqn2.1.1}
\end{equation}

\textit{\noindent Moreover, if the right hand side of (\ref{eqn2.1.1})
exists (finite value), then }$w$\textit{\ }$\in $\textit{\ }$BV_{\vec{F}%
}(\Omega ).$

\bigskip

%TCIMACRO{\TeXButton{Proof}{\proof} }%
%BeginExpansion
\proof
%EndExpansion
For $\vec{\phi}$ $\in $ $C_{0}^{1}(\Omega ),|\vec{\phi}|\leq 1,$ we have%
\begin{eqnarray*}
\int_{\Omega }(-w\func{div}\vec{\phi}+\vec{F}\cdot \vec{\phi})d^{m}x
&=&\lim_{j\rightarrow \infty }\int_{\Omega }(-w_{j}\func{div}\vec{\phi}+\vec{%
F}\cdot \vec{\phi})d^{m}x \\
&\leq &\lim \inf_{j\rightarrow \infty }\int_{\Omega }|Dw_{j}+\vec{F}d^{m}x|
\end{eqnarray*}

\noindent Taking the supremum over all such $\vec{\phi},$ we obtain (\ref%
{eqn2.1.1}) by (\ref{eqn1.2}) (cf. Theorem 5.2.1 in \cite{Zi}). If $\lim
\inf_{j\rightarrow \infty }\int_{\Omega }|Dw_{j}+\vec{F}d^{m}x|$ $<$ $\infty
,$ then $w$\textit{\ }$\in $\textit{\ }$BV_{\vec{F}}(\Omega )$ by definition.

%TCIMACRO{\TeXButton{End Proof}{\endproof}}%
%BeginExpansion
\endproof%
%EndExpansion

Next we claim that $u$ $\in $ $BV(\Omega )$. Since $u_{j}^{-}$ converges to $%
u$ in $C^{0}$-norm (hence $L^{1}$-norm) on $\bar{\Omega},$ we have

\begin{equation}
\int_{\Omega }|Du+\vec{F}d^{m}x|\text{ }\leq \lim \inf_{j\rightarrow \infty
}\int_{\Omega }|\nabla u_{j}^{-}+\vec{F}|d^{m}x\text{.}  \label{eqn2.2}
\end{equation}

\noindent by (\ref{eqn2.1.1}) in Lemma 2.1. We will prove that the right
hand side of (\ref{eqn2.2}) exists (finite value). Let $u_{j,a}^{-}$ denote
the solution of the following elliptic approximating equation:

\begin{eqnarray}
\func{div}(\frac{\nabla \upsilon +\vec{F}}{\sqrt{a^{2}+|\nabla \upsilon +%
\vec{F}|^{2}}}) &=&0\text{ \ in }\Omega  \label{eqn2.3} \\
\upsilon &=&\varphi _{j}^{-}\text{ \ on }\partial \Omega  \notag
\end{eqnarray}

\noindent (cf. (4.1) in \cite{chy}; note that $u_{j,a}^{-}$ $\in $ $%
C^{2,\alpha }$ by Theorem 4.5 in \cite{chy}). From Lemma 2.1 and noting that 
$u_{j}^{-}$ $=$ $\lim_{a\rightarrow 0}u_{j,a}^{-}$ in $C^{0}$-norm (hence $%
L^{1}$-norm) on $\bar{\Omega}$, we have

\begin{equation}
\int_{\Omega }|\nabla u_{j}^{-}+\vec{F}|d^{m}x\leq \lim \inf_{a\rightarrow
0}\int_{\Omega }|\nabla u_{j,a}^{-}+\vec{F}|d^{m}x.  \label{eqn2.4}
\end{equation}

\noindent On the other hand, we observe that

\begin{eqnarray}
|\nabla u_{j,a}^{-}+\vec{F}| &\leq &\sqrt{a^{2}+|\nabla u_{j,a}^{-}+\vec{F}%
|^{2}}  \label{eqn2.5} \\
&=&\nabla u_{j,a}^{-}\cdot N_{j,a}^{-}+\vec{F}\cdot N_{j,a}^{-}+\frac{a^{2}}{%
\sqrt{a^{2}+|\nabla u_{j,a}^{-}+\vec{F}|^{2}}}  \notag
\end{eqnarray}

\noindent where 
\begin{equation}
N_{j,a}^{-}:=\frac{\nabla u_{j,a}^{-}+\vec{F}}{\sqrt{a^{2}+|\nabla
u_{j,a}^{-}+\vec{F}|^{2}}}.  \label{eqn2.6}
\end{equation}

\noindent Integrating (\ref{eqn2.5}) and making use of (\ref{eqn2.3}) (to
get $\nabla u_{j,a}^{-}\cdot N_{j,a}^{-}$ $=$ $\func{div}%
(u_{j,a}^{-}N_{j,a}^{-})$ $-$ $u_{j,a}^{-}\func{div}N_{j,a}^{-}$ $=$ $\func{%
div}(u_{j,a}^{-}N_{j,a}^{-}))$, we obtain

\begin{equation}
\int_{\Omega }|\nabla u_{j,a}^{-}+\vec{F}|d^{m}x\leq \int_{\partial \Omega
}|\varphi _{j}^{-}|d\sigma \text{ }+\text{ }\int_{\Omega }\{|\vec{F}|\text{ }%
+\text{ \TEXTsymbol{\vert}}a|\}d^{m}x  \label{eqn2.7}
\end{equation}

\noindent by noting that $|N_{j,a}^{-}|$ $\leq $ $1,$ where $d\sigma $
denotes the boundary measure$.$ From (\ref{eqn2.7}) we have deduced the
following estimate

\begin{equation}
\lim \inf_{a\rightarrow 0}\int_{\Omega }|\nabla u_{j,a}^{-}+\vec{F}%
|d^{m}x\leq \Vert \varphi _{j}^{-}\Vert _{L^{\infty }(\partial \Omega
)}|\partial \Omega |+\Vert \vec{F}\Vert _{L^{\infty }(\Omega )}|\Omega |
\label{eqn2.8}
\end{equation}

\noindent where $|\partial \Omega |$ and $|\Omega |$ denote the $2n-1$ and $%
2n$ dimensional Hausdorff measures of $\partial \Omega $ and $\Omega ,$
respectively. It now follows from (\ref{eqn2.2}), (\ref{eqn2.4}), and (\ref%
{eqn2.8}) that

\begin{equation}
\int_{\Omega }|Du+\vec{F}d^{m}x|\leq \Vert \varphi \Vert _{L^{\infty
}(\partial \Omega )}|\partial \Omega |+\Vert \vec{F}\Vert _{L^{\infty
}(\Omega )}|\Omega |<\infty .  \label{eqn2.9}
\end{equation}

\noindent (\ref{eqn2.9}) means $u$ $\in $ $BV(\Omega ).$ In the remaining
section, we will show that $u$ is a minimizer for $\mathcal{F}(\cdot )$ in $%
C^{0}(\bar{\Omega})\cap BV(\Omega )$ with the same boundary value $\varphi $%
. Take $v$ $\in $ $C^{0}(\bar{\Omega})\cap BV(\Omega )$ such that $v$ $=$ $u$
on $\partial \Omega .$ Let $v_{\tau }$ be a mollifier of $v,$ where $\tau $ $%
>$ $0$. Let $\vec{G}_{\tau }$ $=$ $((G_{1})_{\tau },$ $(G_{2})_{\tau },$ $%
...)$ be a mollifier of a vector field $\vec{G}$ $=$ $(G_{1},$ $G_{2},$ $%
...).$

\bigskip

\textbf{Lemma 2.2}. \textit{Let }$\Omega ^{\prime }$\textit{\ }$\subset
\subset $\textit{\ }$\Omega $\textit{\ (i.e., }$\Omega ^{\prime }$\textit{\
has compact closure in }$\Omega ).$\textit{\ For }$\tau $\textit{\ small
enough, there holds}

\begin{equation}
\int_{\Omega ^{\prime }}|\nabla v_{\tau }+\vec{F}|d^{m}x\leq \int_{\Omega
}|Dv+\vec{F}d^{m}x|\text{ }+\text{ }\Vert \vec{F}-\vec{F}_{\tau }\Vert
_{L^{1}(\Omega )}\mathit{.}  \label{eqn2.10}
\end{equation}

\bigskip

%TCIMACRO{\TeXButton{Proof}{\proof} }%
%BeginExpansion
\proof
%EndExpansion
For $\vec{\phi}$ $\in $ $C_{0}^{1}(\Omega ^{\prime })$ such that $|\vec{\phi}%
|$ $\leq $ $1$ (which implies $|\vec{\phi}_{\tau }|$ $\leq $ $1),$ we compute%
\begin{eqnarray}
&&\int_{\Omega ^{\prime }}(-v_{\tau }\func{div}\vec{\phi}+\vec{F}\cdot \vec{%
\phi})d^{m}x  \label{eqn2.11} \\
&=&\int_{\Omega }(-v\func{div}\vec{\phi}_{\tau }+\vec{F}\cdot \vec{\phi}%
_{\tau }+\vec{F}\cdot (\vec{\phi}-\vec{\phi}_{\tau }))d^{m}x  \notag \\
&\leq &\int_{\Omega }|Dv+\vec{F}d^{m}x|\text{ }+\text{ }\int_{\Omega }(\vec{F%
}-\vec{F}_{\tau })\cdot \vec{\phi}\text{ }d^{m}x  \notag \\
&\leq &\int_{\Omega }|Dv+\vec{F}d^{m}x|\text{ }+\text{ }\Vert \vec{F}-\vec{F}%
_{\tau }\Vert _{L^{1}(\Omega )}.  \notag
\end{eqnarray}

\noindent By taking the supremum of the left side of (\ref{eqn2.11}) over $%
\vec{\phi}$, we obtain (\ref{eqn2.10}).

%TCIMACRO{\TeXButton{End Proof}{\endproof}}%
%BeginExpansion
\endproof%
%EndExpansion

\bigskip

\textbf{Lemma 2.3.} \textit{Let }$\upsilon ,$\textit{\ }$\omega $\textit{\ }$%
\in $\textit{\ }$C^{2}(\bar{\Omega})$\textit{\ satisfy }$\func{div}%
N_{a}(\upsilon )=\func{div}N_{a}(\omega )=0$\textit{\ in }$\Omega ,$ \textit{%
where} $N_{a}(\rho )$ $:=$ $\frac{\nabla \rho +\vec{F}}{\sqrt{a^{2}+|\nabla
\rho +\vec{F}|^{2}}}$\textit{. For any }$\Omega ^{\prime }$\textit{\ }$%
\subset \subset $\textit{\ }$\Omega $\textit{, there holds}%
\begin{equation}
\mid \int_{\Omega ^{\prime }}\{\sqrt{a^{2}+|\nabla \upsilon +\vec{F}|^{2}}-%
\sqrt{a^{2}+|\nabla \omega +\vec{F}|^{2}}\}d^{m}x\mid \leq \int_{\partial
\Omega ^{\prime }}|\upsilon -\omega |d\sigma .  \label{eqn2.12}
\end{equation}

\bigskip

%TCIMACRO{\TeXButton{Proof}{\proof} }%
%BeginExpansion
\proof
%EndExpansion
Consider the following expression

\begin{equation}
I(s):=\int_{\Omega ^{\prime }}\sqrt{a^{2}+|\nabla \upsilon +\vec{F}+s\nabla
(\omega -\upsilon )|^{2}}d^{m}x+s\int_{\partial \Omega ^{\prime }}|\upsilon
-\omega |d\sigma   \label{eqn2.13}
\end{equation}

\noindent for $s$ $\in $ $[0,1].$ Compute

\begin{equation}
I^{\prime }(s)=\int_{\Omega ^{\prime }}\frac{[\nabla \upsilon +\vec{F}%
+s\nabla (\omega -\upsilon )]\cdot \nabla (\omega -\upsilon )}{\sqrt{%
a^{2}+|\nabla \upsilon +\vec{F}+s\nabla (\omega -\upsilon )|^{2}}}%
d^{m}x+\int_{\partial \Omega ^{\prime }}|\upsilon -\omega |d\sigma
\label{eqn2.14}
\end{equation}

\noindent and 
\begin{eqnarray*}
&&I^{\prime \prime }(s) \\
&=&\int_{\Omega ^{\prime }}\{\frac{|\nabla (\omega -\upsilon )|^{2}}{\sqrt{%
a^{2}+|\nabla \upsilon +\vec{F}+s\nabla (\omega -\upsilon )|^{2}}}- \\
&&\frac{\{[\nabla \upsilon +\vec{F}+s\nabla (\omega -\upsilon )]\cdot \nabla
(\omega -\upsilon )\}^{2}}{(\sqrt{a^{2}+|\nabla \upsilon +\vec{F}+s\nabla
(\omega -\upsilon )|^{2}})^{3}}\}d^{m}x \\
&=&\int_{\Omega ^{\prime }}\{a^{2}|\nabla (\omega -\upsilon )|^{2}+|\nabla
(\omega -\upsilon )|^{2}|\nabla \upsilon +\vec{F}+s\nabla (\omega -\upsilon
)|^{2} \\
&&-\{[\nabla \upsilon +\vec{F}+s\nabla (\omega -\upsilon )]\cdot \nabla
(\omega -\upsilon )\}^{2}\}(\sqrt{a^{2}+|\nabla \upsilon +\vec{F}+s\nabla
(\omega -\upsilon )|^{2}})^{-3}d^{m}x \\
&\geq &0
\end{eqnarray*}

\noindent by Cauchy's inequality for $s$ $\in $ $[0,1]$. It follows that

\begin{equation}
I^{\prime }(s)\geq I^{\prime }(0)  \label{eqn2.15}
\end{equation}

On the other hand, from (\ref{eqn2.14}) we compute%
\begin{eqnarray}
I^{\prime }(0) &=&\int_{\Omega ^{\prime }}N_{a}(\upsilon )\cdot \nabla
(\omega -\upsilon )d^{m}x+\int_{\partial \Omega ^{\prime }}|\upsilon -\omega
|d\sigma  \label{eqn2.16} \\
&=&\int_{\partial \Omega ^{\prime }}(\omega -\upsilon )N_{a}(\upsilon )\cdot
\nu d\sigma +\int_{\partial \Omega ^{\prime }}|\upsilon -\omega |d\sigma 
\notag \\
&\geq &0  \notag
\end{eqnarray}

\noindent where we have used the equation $\func{div}N_{a}(\upsilon )$ $=$ $%
0 $ and $|N_{a}(\upsilon )|$ $\leq $ $1.$ By (\ref{eqn2.15}) and (\ref%
{eqn2.16}), we get $I^{\prime }(s)$ $\geq $ $0,$ and hence $I(1)\geq I(0).$
That is%
\begin{equation*}
\int_{\Omega ^{\prime }}\{\sqrt{a^{2}+|\nabla \upsilon +\vec{F}|^{2}}-\sqrt{%
a^{2}+|\nabla \omega +\vec{F}|^{2}}\}d^{m}x\leq \int_{\partial \Omega
^{\prime }}|\upsilon -\omega |d\sigma .
\end{equation*}

\noindent Switching $\upsilon $ and $\omega $ in the above argument, we
finally reach (\ref{eqn2.12}).

%TCIMACRO{\TeXButton{End Proof}{\endproof}}%
%BeginExpansion
\endproof%
%EndExpansion

\bigskip

%TCIMACRO{\TeXButton{Proof}{\proof} }%
%BeginExpansion
\proof
%EndExpansion
(\textbf{of Theorem A continued)} Now we consider only parabolically convex
domain $\Omega ^{\prime }$ $\subset \subset $ $\Omega $ with $\partial
\Omega ^{\prime }$ $\in $ $C^{\infty }.$ For $a>0$ let $u_{j,a}^{-},$ $%
v_{\tau ,a}$ be the solutions to $\func{div}N_{a}(\cdot )$ $=$ $0$ in $%
\Omega ^{\prime },$ such that $u_{j,a}^{-}$ $=$ $u_{j}^{-},$ $v_{\tau ,a}$ $%
= $ $v_{\tau }$ on $\partial \Omega ^{\prime }$, where $N_{a}(\rho )$ $:=$ $%
\frac{\nabla \rho +\vec{F}}{\sqrt{a^{2}+|\nabla \rho +\vec{F}|^{2}}}.$ We
then compute%
\begin{eqnarray}
&&\int_{\Omega ^{\prime }}|\nabla u_{j}^{-}+\vec{F}|d^{m}x  \label{eqn2.17}
\\
&\leq &\int_{\Omega ^{\prime }}|\nabla u_{j,a}^{-}+\vec{F}|d^{m}x\text{ \ (}%
u_{j}^{-}\text{ is a minimizer for }\mathcal{F}(\cdot ))  \notag \\
&\leq &\int_{\Omega ^{\prime }}\sqrt{a^{2}+|\nabla u_{j,a}^{-}+\vec{F}|^{2}}%
d^{m}x  \notag \\
&\leq &\int_{\Omega ^{\prime }}\sqrt{a^{2}+|\nabla v_{\tau ,a}+\vec{F}|^{2}}%
d^{m}x+\int_{\partial \Omega ^{\prime }}|u_{j}^{-}-v_{\tau }|d\sigma  \notag
\end{eqnarray}

\noindent by (\ref{eqn2.12}). Let $\mathcal{F}_{a}(w)$ $\equiv $ $%
\int_{\Omega ^{\prime }}\sqrt{a^{2}+|\nabla w+\vec{F}|^{2}}d^{m}x.$ Since $%
v_{\tau ,a}$ is a minimizer for $\mathcal{F}_{a}(\cdot )$ (see \cite{chy}),
we estimate%
\begin{eqnarray}
&&\int_{\Omega ^{\prime }}\sqrt{a^{2}+|\nabla v_{\tau ,a}+\vec{F}|^{2}}d^{m}x
\label{eqn2.18} \\
&\leq &\int_{\Omega ^{\prime }}\sqrt{a^{2}+|\nabla v_{\tau }+\vec{F}|^{2}}%
d^{m}x  \notag \\
&\leq &a\text{ }|\Omega ^{\prime }|\text{ }+\text{ }\int_{\Omega }|\nabla v+%
\vec{F}|d^{m}x\text{ }+\text{ }\Vert \vec{F}-\vec{F}_{\tau }\Vert
_{L^{1}(\Omega )}  \notag
\end{eqnarray}

\noindent by (\ref{eqn2.10}) ($\tau $ small enough). Combining (\ref{eqn2.17}%
) and (\ref{eqn2.18}) gives%
\begin{eqnarray}
&&\int_{\Omega ^{\prime }}|\nabla u_{j}^{-}+\vec{F}|d^{m}x  \label{eqn2.19}
\\
&\leq &a\text{ }|\Omega ^{\prime }|\text{ }+\text{ }\int_{\Omega }|\nabla v+%
\vec{F}|d^{m}x\text{ }+\text{ }\Vert \vec{F}-\vec{F}_{\tau }\Vert
_{L^{1}(\Omega )}+\int_{\partial \Omega ^{\prime }}|u_{j}^{-}-v_{\tau
}|d\sigma .  \notag
\end{eqnarray}

\noindent In view of (\ref{eqn2.2}) and (\ref{eqn2.19}), we conclude (\ref%
{eqn1.4}) by letting $a$ go to zero and $\Omega ^{\prime }$ approach $\Omega 
$ ($\tau $ tends to zero accordingly).

%TCIMACRO{\TeXButton{End Proof}{\endproof}}%
%BeginExpansion
\endproof%
%EndExpansion

\bigskip

\section{\textbf{Extension to measures and the first variation}}

We will extend the domain of $\mathcal{F(\cdot )}$ from $BV$ functions to
vector-valued measures. Let $\Omega $ $\subset $ $R^{m}$ be a bounded
domain. Let $u$ $\in $ $BV_{\vec{F}}(\Omega )$ where $\vec{F}$ $:=$ $(F_{i})$
$\in $ $L^{1}(\Omega ).$ Then there are Radon signed measures $\lambda _{1},$
$\lambda _{2},$ ...,$\lambda _{m}$ defined in $\Omega $ such that for $i$ $=$
$1,$ $2,$...$,m,$ (Recall that $d^{m}x$ denotes the Lebesgue measure of $%
R^{m})$%
\begin{equation*}
\int_{\Omega }u\frac{\partial \varphi }{\partial x_{i}}d^{m}x=-\int_{\Omega
}\varphi d\lambda _{i}
\end{equation*}

\noindent for all $\varphi $ $\in $ $C_{0}^{\infty }(\Omega )$ (see Remark
1.5 on page 5 in \cite{Gi} or see (5.1.1) in \cite{Zi}, and note that $u$ $%
\in $ $BV_{\vec{F}}(\Omega )$ if and only if $u$ $\in $ $BV(\Omega )$)$.$
Write $Du$ $:=$ $(d\lambda _{i}).$ So $Du+\vec{F}d^{m}x$ defines a
vector-valued Radon signed measure and we define its total variation
(measure)

\begin{equation*}
|Du+\vec{F}d^{m}x|(f):=\sup \{\int_{\Omega }(-u\func{div}\vec{\phi}+\vec{F}%
\cdot \vec{\phi})d^{m}x\mid \vec{\phi}\in C_{0}^{1}(\Omega ),|\vec{\phi}%
|\leq f\}
\end{equation*}

\noindent for $f$ being a non-negative real-valued continuous function with
compact support in $\Omega .$ By the Riesz Representation Theorem, $|Du+\vec{%
F}d^{m}x|$ is a non-negative Radon measure on $\Omega $ (mimicking the
argument in Remark 5.1.2. of \cite{Zi})$.$ Similarly, for a general
vector-valued measure $d\mu $ $=$ $(d\mu _{i})$ (instead of $\mu $ $=$ $(\mu
_{i}))$ on $\Omega ,$ we define its total variation measure $|d\mu |$ by 
\begin{equation*}
|d\mu |(f)=\sup \{\int_{\Omega }\vec{\phi}\cdot d\mu \mid \vec{\phi}\in
C_{0}^{1}(\Omega ),|\vec{\phi}|\leq f\}.
\end{equation*}

\noindent We extend the domain of $\mathcal{F(\cdot )}$ to include
vector-valued (Radon signed) measures $d\mu $ by defining

\begin{equation}
\mathcal{F(}d\mu ):=\int_{\Omega }|d\mu |:=|d\mu |(1).  \label{eqn3.1}
\end{equation}

\noindent In this section we want to compute the first variation of $%
\mathcal{F(\cdot )}$ in measures.

Let $E$ be a $C^{\infty }$-smooth Riemannian vector bundle over a $C^{\infty
}$-smooth manifold $X$. Let $d\mu ,$ $d\nu $ be two $E$-valued measures on $%
X $. We assume that both $d\mu $ and $d\nu $ are bounded in the sense that $%
|d\mu |$ and $|d\nu |$ are integrable over $X,$ i.e., $\mathcal{F(}d\mu )$
and $\mathcal{F(}d\nu )$ are finite in view of (\ref{eqn3.1}) with $\Omega $
replaced by $X$ ($\vec{\phi}$ is viewed as a $C^{1}$-smooth section of $E$
with compact support while "$\cdot "$ denotes the fibre inner product). Let $%
d\mu _{\varepsilon }$ $:=$ $d\mu $ $+$ $\varepsilon d\nu $ for $\varepsilon $
$\in $ $R$. Since $|d\mu _{\varepsilon }|$ is a positive bounded measure, we
can find $N_{\varepsilon }$, $A_{\varepsilon }$ $\in $ $L^{1}(|d\mu
_{\varepsilon }|)$ with $|N_{\varepsilon }|$ $=$ $1,$ such that

\begin{eqnarray}
d\mu _{\varepsilon } &=&N_{\varepsilon }|d\mu _{\varepsilon }|,
\label{eqn3.2} \\
d\nu &=&A_{\varepsilon }|d\mu _{\varepsilon }|+d\nu _{s}^{\varepsilon },%
\text{ }d\nu _{s}^{\varepsilon }\perp |d\mu _{\varepsilon }|  \notag
\end{eqnarray}

\noindent according to the Radon-Nikodym theorem (extending 6.9 and 6.12 in 
\cite{Ru} to the case of vector-valued measures; see also \cite{Ro}).

\bigskip

%TCIMACRO{\TeXButton{Proof}{\proof} }%
%BeginExpansion
\proof
%EndExpansion
\textbf{(of Theorem B)} We have

\begin{eqnarray}
|d\mu _{\varepsilon _{2}}| &=&|d\mu _{\varepsilon _{1}}+(\varepsilon
_{2}-\varepsilon _{1})d\nu |  \label{eqn3.3} \\
&=&|(\varepsilon _{2}-\varepsilon _{1})A_{\varepsilon _{1}}+N_{\varepsilon
_{1}}||d\mu _{\varepsilon _{1}}|+|\varepsilon _{2}-\varepsilon _{1}||d\nu
_{s}^{\varepsilon _{1}}|  \notag
\end{eqnarray}

\noindent by (\ref{eqn3.2}). It follows from (\ref{eqn3.3}) that for $%
\varepsilon _{2}\neq \varepsilon _{1}$

\begin{eqnarray}
\frac{|d\mu _{\varepsilon _{2}}|-|d\mu _{\varepsilon _{1}}|}{\varepsilon
_{2}-\varepsilon _{1}} &=&\{|(\varepsilon _{2}-\varepsilon
_{1})A_{\varepsilon _{1}}+N_{\varepsilon _{1}}|-1\}\frac{|d\mu _{\varepsilon
_{1}}|}{\varepsilon _{2}-\varepsilon _{1}}  \label{eqn3.4} \\
&&+\frac{|\varepsilon _{2}-\varepsilon _{1}|}{\varepsilon _{2}-\varepsilon
_{1}}|d\nu _{s}^{\varepsilon _{1}}|.  \notag
\end{eqnarray}

\noindent Observe that 
\begin{eqnarray}
&&|\frac{|(\varepsilon _{2}-\varepsilon _{1})A_{\varepsilon
_{1}}+N_{\varepsilon _{1}}|-1}{\varepsilon _{2}-\varepsilon _{1}}||d\mu
_{\varepsilon _{1}}|  \label{3.4.1} \\
&=&|\frac{|(\varepsilon _{2}-\varepsilon _{1})A_{\varepsilon
_{1}}+N_{\varepsilon _{1}}|-|N_{\varepsilon _{1}}|}{\varepsilon
_{2}-\varepsilon _{1}}||d\mu _{\varepsilon _{1}}|  \notag \\
&\leq &|A_{\varepsilon _{1}}||d\mu _{\varepsilon _{1}}|\leq |d\nu |  \notag
\end{eqnarray}%
\noindent by noting that $|N_{\varepsilon _{1}}|$ $=$ $1.$ It follows from (%
\ref{eqn3.4}) and (\ref{3.4.1}) that $\mathcal{F(\varepsilon )}$ is
Lipschitz continuous in $\varepsilon $\textit{\ }since $d\nu $ is bounded by
assumption. Also observe that

\begin{eqnarray}
|(\varepsilon _{2}-\varepsilon _{1})A_{\varepsilon _{1}}+N_{\varepsilon
_{1}}|-1 &=&\frac{|(\varepsilon _{2}-\varepsilon _{1})A_{\varepsilon
_{1}}+N_{\varepsilon _{1}}|^{2}-1}{|(\varepsilon _{2}-\varepsilon
_{1})A_{\varepsilon _{1}}+N_{\varepsilon _{1}}|+1}  \label{eqn3.5} \\
&=&\frac{2(\varepsilon _{2}-\varepsilon _{1})A_{\varepsilon _{1}}\cdot
N_{\varepsilon _{1}}+|(\varepsilon _{2}-\varepsilon _{1})A_{\varepsilon
_{1}}|^{2}}{|(\varepsilon _{2}-\varepsilon _{1})A_{\varepsilon
_{1}}+N_{\varepsilon _{1}}|+1}.  \notag
\end{eqnarray}

\noindent Since $|A_{\varepsilon _{1}}||d\mu _{\varepsilon _{1}}|$ and $%
|d\nu _{s}^{\varepsilon _{1}}|$ are integrable by assumption ($d\nu $ is
bounded), we can invoke the Lebesgue dominated convergence theorem to obtain

\begin{equation}
\lim_{\varepsilon _{2}\rightarrow \varepsilon _{1}\pm }\int_{X}\frac{|d\mu
_{\varepsilon _{2}}|-|d\mu _{\varepsilon _{1}}|}{\varepsilon
_{2}-\varepsilon _{1}}=\int_{X}N_{\varepsilon _{1}}\cdot A_{\varepsilon
_{1}}|d\mu _{\varepsilon _{1}}|\pm \int_{X}|d\nu _{s}^{\varepsilon _{1}}|
\label{eqn3.6}
\end{equation}

\noindent by (\ref{eqn3.4}) and (\ref{eqn3.5}).

%TCIMACRO{\TeXButton{End Proof}{\endproof}}%
%BeginExpansion
\endproof%
%EndExpansion

\bigskip

%TCIMACRO{\TeXButton{Proof}{\proof} }%
%BeginExpansion
\proof
%EndExpansion
\textbf{(of Corollary B$^{\prime}$)} Let $d\mu $ $=$ $Du$ $+$ $\vec{F}d^{m}x$
denote the vector-valued measure associated to $u$ $\in $ $BV(\Omega )$. Let 
$d\nu $ $=$ $D\varphi $ for $\varphi $ $\in $ $BV(\Omega )$ with $\varphi
|_{\partial \Omega }$ $=$ $0.$ Recall that we denote $\mathcal{\tilde{F}}%
_{H}(u+\varepsilon \varphi )$ by $\mathcal{\tilde{F}}_{H}(\varepsilon )$.%
\textit{\ }Now\textit{\ }it is straightforward to extend (\ref{eqn1.5}) for $%
X$ $=$ $\Omega $ to include $H$ as below:%
\begin{equation}
\mathcal{\tilde{F}}_{H}^{\prime }(\mathcal{\varepsilon }_{1}\mathcal{\pm }%
)=\int_{\Omega }N_{\varepsilon _{1}}\cdot A_{\varepsilon _{1}}|d\mu
_{\varepsilon _{1}}|\pm \int_{\Omega }|d\nu _{s}^{\varepsilon
_{1}}|+\int_{\Omega }H\varphi \text{ }d^{m}x.  \label{eqn3.6.1}
\end{equation}%
\noindent Letting $\varepsilon _{1}$ $=$ $0$ in (\ref{eqn3.6.1}) we have 
\begin{eqnarray*}
&&\int_{\Omega }N_{0}\cdot A_{0}|d\mu |\pm \int_{\Omega }|d\nu
_{s}^{0}|+\int_{\Omega }H\varphi \text{ }d^{m}x \\
&=&\mathcal{\tilde{F}}_{H}^{\prime }(0\mathcal{\pm }) \\
&=&\lim_{\varepsilon \rightarrow 0\pm }\frac{\mathcal{\tilde{F}}%
_{H}(u+\varepsilon \varphi )-\mathcal{\tilde{F}}_{H}(u)}{\varepsilon }\geq 0%
\text{ (}\leq 0,\text{ resp.})
\end{eqnarray*}

\noindent for $\varepsilon $ $\rightarrow $ $0+$ ($\varepsilon $ $%
\rightarrow $ $0-,$ resp.) since $\mathcal{\tilde{F}}_{H}(u+\varepsilon
\varphi )$ $-$ $\mathcal{\tilde{F}}_{H}(u)\geq 0$ for $u$ being a minimizer
and $\varepsilon $ $>$ $0$ ($\varepsilon $ $<$ $0,$ resp.)$.$ We have proved
(\ref{eqn1.5.0}) and (\ref{eqn1.5.1}).

%TCIMACRO{\TeXButton{End Proof}{\endproof}}%
%BeginExpansion
\endproof%
%EndExpansion

\bigskip

\textbf{Lemma 3.1}. \textit{Suppose }$d\mu ,$\textit{\ }$d\nu $\textit{\ are
two bounded }$E$\textit{-valued measures on }$X$\textit{\ as described
between (\ref{eqn3.1}) and (\ref{eqn3.2}). Let }$d\mu _{\varepsilon }$%
\textit{\ }$:=$\textit{\ }$d\mu $\textit{\ }$+$\textit{\ }$\varepsilon d\nu $%
\textit{\ for }$\varepsilon $\textit{\ }$\in $\textit{\ }$R$\textit{\
satisfy (\ref{eqn3.2}). Then for }$\varepsilon _{1}$\textit{\ }$\neq $%
\textit{\ }$\varepsilon _{2}$\textit{\ there holds }$d\nu _{s}^{\varepsilon
_{1}}\perp d\nu _{s}^{\varepsilon _{2}},$\textit{\ i.e., }$\mid d\nu
_{s}^{\varepsilon _{1}}\mid \perp \mid d\nu _{s}^{\varepsilon _{2}}\mid .$%
\textit{\ Moreover, there exist at most countably many }$\varepsilon $%
\textit{'s such that \TEXTsymbol{\vert}}$d\nu _{s}^{\varepsilon }|(X)$%
\textit{\ }$\neq $\textit{\ }$0.$

\bigskip

%TCIMACRO{\TeXButton{Proof}{\proof} }%
%BeginExpansion
\proof
%EndExpansion
Let $j$ $=$ $1$, $2.$ Since $|d\mu _{\varepsilon _{j}}|$ $\perp $ $|d\nu
_{s}^{\varepsilon _{j}}|,$ we can find a measurable set $E_{\varepsilon
_{j}} $ such that $|d\mu _{\varepsilon _{j}}|$ is concentrated on $%
E_{\varepsilon _{j}}$ and $|d\nu _{s}^{\varepsilon _{j}}|$ is concentrated
on $E_{\varepsilon _{j}}^{c},$ the complement of $E_{\varepsilon _{j}}.$ For
any measurable set $B$ $\subset $ $E_{\varepsilon _{1}}^{c}\cap
E_{\varepsilon _{2}}^{c},$ \TEXTsymbol{\vert}$d\nu |(B)$ $=$ $0$ by
observing that $d\mu _{\varepsilon _{1}}$ $-$ $d\mu _{\varepsilon _{2}}$ $=$ 
$(\varepsilon _{1}-\varepsilon _{2})d\nu .$ It follows that \TEXTsymbol{\vert%
}$d\nu _{s}^{\varepsilon _{j}}|(B)$ $=$ $0$ for $j$ $=$ $1,$ $2$ since $d\nu 
$ $=$ $d\nu _{s}^{\varepsilon _{1}}$ and $d\nu $ $=$ $d\nu _{s}^{\varepsilon
_{2}}$ on $E_{\varepsilon _{1}}^{c}\cap E_{\varepsilon _{2}}^{c}.$ So $|d\nu
_{s}^{\varepsilon _{1}}|$ $=$ $|d\nu _{s}^{\varepsilon _{2}}|$ $=$ $0$ on $%
E_{\varepsilon _{1}}^{c}\cap E_{\varepsilon _{2}}^{c},$ and hence $|d\nu
_{s}^{\varepsilon _{1}}|$ and $|d\nu _{s}^{\varepsilon _{2}}|$ are
concentrated on $E_{\varepsilon _{1}}^{c}\backslash $ ($E_{\varepsilon
_{1}}^{c}\cap E_{\varepsilon _{2}}^{c})$ and $E_{\varepsilon
_{2}}^{c}\backslash $ ($E_{\varepsilon _{1}}^{c}\cap E_{\varepsilon
_{2}}^{c})$ (the intersection of these two sets is empty), respectively.
Therefore $|d\nu _{s}^{\varepsilon _{1}}|$ $\perp $ $|d\nu _{s}^{\varepsilon
_{2}}|\mathit{.}$

Given a positive integer $n,$ we can only have finitely many $\varepsilon
_{j}$'s such that $|d\nu |(E_{\varepsilon _{j}}^{c})$ $=$ $|d\nu
_{s}^{\varepsilon _{j}}|(E_{\varepsilon _{j}}^{c})\geq \frac{1}{n}$ since 
\begin{eqnarray*}
\sum_{j}|d\nu |(E_{\varepsilon _{j}}^{c}) &=&\sum_{j}|d\nu _{s}^{\varepsilon
_{j}}|(E_{\varepsilon _{j}}^{c}) \\
&=&|d\nu |(\cup E_{\varepsilon _{j}}^{c}) \\
&\leq &|d\nu |(X)<\infty \text{ \ (by assumption).}
\end{eqnarray*}

\noindent It follows that there are at most countably many $\varepsilon $'s
such that $|d\nu |(E_{\varepsilon }^{c})$ $=$ $|d\nu _{s}^{\varepsilon
}|(E_{\varepsilon }^{c})$ $\neq $ $0.$

%TCIMACRO{\TeXButton{End Proof}{\endproof}}%
%BeginExpansion
\endproof%
%EndExpansion

\bigskip

If $\varepsilon $ satisfies \textit{\TEXTsymbol{\vert}}$d\nu
_{s}^{\varepsilon }|(X)$\textit{\ }$\neq $\textit{\ }$0,$ we call it
singular, otherwise regular. Denote $\mathcal{F(}d\mu _{\varepsilon })$ by $%
\mathcal{F(\varepsilon )}$ for simplicity.

\bigskip

\textbf{Lemma 3.2.} \textit{Suppose we are in the situation of Lemma 3.1.
Then we have}

\textit{(1) For }$\varepsilon _{1},$\textit{\ }$\varepsilon _{2}$\textit{\
arbitrary, there holds}%
\begin{equation}
|d\mu _{\varepsilon _{1}}|+|d\nu _{s}^{\varepsilon _{1}}|\ll |d\mu
_{\varepsilon _{2}}|+|d\nu _{s}^{\varepsilon _{2}}|\ll |d\mu _{\varepsilon
_{1}}|+|d\nu _{s}^{\varepsilon _{1}}|;  \label{3.7}
\end{equation}

\textit{(2) For }$\varepsilon _{1},$\textit{\ }$\varepsilon _{2}$\textit{\
regular, there holds }$|d\mu _{\varepsilon _{1}}|$\textit{\ }$\ll $\textit{\ 
}$|d\mu _{\varepsilon _{2}}|$\textit{\ }$\ll $\textit{\ }$|d\mu
_{\varepsilon _{1}}|.$

\bigskip

%TCIMACRO{\TeXButton{Proof}{\proof} }%
%BeginExpansion
\proof
%EndExpansion
We may assume $\varepsilon _{1}$ $\neq $\textit{\ }$\varepsilon _{2}.$ From (%
\ref{eqn3.3}) we have

\begin{equation}
|d\mu _{\varepsilon _{2}}|\ll |d\mu _{\varepsilon _{1}}|+|d\nu
_{s}^{\varepsilon _{1}}|  \label{3.8}
\end{equation}

\noindent Switching $\varepsilon _{1}$ and\textit{\ }$\varepsilon _{2}$ in (%
\ref{eqn3.3}) gives%
\begin{eqnarray}
&&|d\mu _{\varepsilon _{1}}|  \label{3.9} \\
&=&|(\varepsilon _{1}-\varepsilon _{2})A_{\varepsilon _{2}}+N_{\varepsilon
_{2}}||d\mu _{\varepsilon _{2}}|+|\varepsilon _{1}-\varepsilon _{2}||d\nu
_{s}^{\varepsilon _{2}}|.  \notag
\end{eqnarray}%
\noindent Therefore we obtain%
\begin{equation}
|\varepsilon _{1}-\varepsilon _{2}||d\nu _{s}^{\varepsilon _{2}}|\ll |d\mu
_{\varepsilon _{1}}|+|d\nu _{s}^{\varepsilon _{1}}|  \label{3.10}
\end{equation}%
\noindent by (\ref{3.9}) and (\ref{3.8}). Now it follows from (\ref{3.8}%
).and (\ref{3.10}) that%
\begin{equation*}
|d\mu _{\varepsilon _{2}}|+|d\nu _{s}^{\varepsilon _{2}}|\ll |d\mu
_{\varepsilon _{1}}|+|d\nu _{s}^{\varepsilon _{1}}|.
\end{equation*}%
\noindent By symmetry (\ref{3.7}) follows. For $\varepsilon _{1},$\textit{\ }%
$\varepsilon _{2}$ regular, $|d\nu _{s}^{\varepsilon _{1}}|$ $=$ $|d\nu
_{s}^{\varepsilon _{2}}|$ $=$ $0$ and hence (2) follows from (\ref{3.7}).

%TCIMACRO{\TeXButton{End Proof}{\endproof}}%
%BeginExpansion
\endproof%
%EndExpansion

\bigskip

\textbf{Lemma 3.3}. \textit{Suppose }$d\mu $\textit{\ and }$d\mu ^{\prime }$%
\textit{\ are two bounded }$E$\textit{-valued measures on }$X$\textit{\ (see
the paragraph between (\ref{eqn3.1}) and (\ref{eqn3.2})). Assume \TEXTsymbol{%
\vert}}$d\mu |$\textit{\ }$\ll $\textit{\ \TEXTsymbol{\vert}}$d\mu ^{\prime
}|$\textit{\ }$\ll $\textit{\ \TEXTsymbol{\vert}}$d\mu |.$\textit{\ Write }$%
d\mu $\textit{\ }$=$\textit{\ }$N_{\mu }|d\mu |,$\textit{\ }$d\mu ^{\prime }$%
\textit{\ }$=$\textit{\ }$N_{\mu ^{\prime }}|d\mu ^{\prime }|.$\textit{\
Then we have }

\textit{(1) }$N_{\mu }$\textit{\ }$\neq $\textit{\ }$0$\textit{\ a.e. }[$%
|d\mu ^{\prime }|]$\textit{\ and }$N_{\mu ^{\prime }}$\textit{\ }$\neq $%
\textit{\ }$0$\textit{\ a.e. }[$|d\mu |]$\textit{\ and}

\textit{(2) there holds}%
\begin{equation}
(N_{\mu }-N_{\mu ^{\prime }})\cdot (d\mu -d\mu ^{\prime })=\frac{1}{2}%
|N_{\mu }-N_{\mu ^{\prime }}|^{2}(|d\mu |+|d\mu ^{\prime }|).
\label{eqn3.8.1}
\end{equation}

\bigskip

%TCIMACRO{\TeXButton{Proof}{\proof} }%
%BeginExpansion
\proof
%EndExpansion
Suppose there is a set $S$ with $|d\mu ^{\prime }|(S)$ $>$ $0$ and $N_{\mu }$
$=$ $0$ on $S.$ From the definition of $N_{\mu },$ we have \TEXTsymbol{\vert}%
$d\mu |(S)$ $=$ $0.$ It follows that \TEXTsymbol{\vert}$d\mu ^{\prime }|(S)$ 
$=$ $0$ by the assumption \TEXTsymbol{\vert}$d\mu ^{\prime }|$ $\ll $ 
\TEXTsymbol{\vert}$d\mu |.$ We have reached a contradiction. Therefore $%
N_{\mu }$ $\neq $ $0$ a.e. [$|d\mu ^{\prime }|].$ By symmetry, we also have $%
N_{\mu ^{\prime }}$ $\neq $ $0$ a.e. [$|d\mu |].$ We have proved (1).

As for (2), noting that $N_{\mu },$ $N_{\mu ^{\prime }}$ are defined a.e. [$%
|d\mu |]$ and [$|d\mu ^{\prime }|],$ we compute%
\begin{eqnarray}
(N_{\mu }-N_{\mu ^{\prime }})\cdot (d\mu -d\mu ^{\prime }) &=&(N_{\mu
}-N_{\mu ^{\prime }})(N_{\mu }|d\mu |-N_{\mu ^{\prime }}|d\mu ^{\prime }|)
\label{eqn3.8.2} \\
&=&(1-N_{\mu }\cdot N_{\mu ^{\prime }})(|d\mu |+|d\mu ^{\prime }|)  \notag \\
&=&\frac{1}{2}|N_{\mu }-N_{\mu ^{\prime }}|^{2}(|d\mu |+|d\mu ^{\prime }|). 
\notag
\end{eqnarray}

%TCIMACRO{\TeXButton{End Proof}{\endproof}}%
%BeginExpansion
\endproof%
%EndExpansion

\bigskip

We remark that for general $d\mu ,$ $d\mu ^{\prime }$ (which may not satisfy
the condition \TEXTsymbol{\vert}$d\mu |$\textit{\ }$\ll $\textit{\ 
\TEXTsymbol{\vert}}$d\mu ^{\prime }|$\textit{\ }$\ll $\textit{\ \TEXTsymbol{%
\vert}}$d\mu |),$ the formula (\ref{eqn3.8.1}) should be interpreted and
modified as below. Write $d\mu $ $=$ $A|d\mu ^{\prime }|$ $+$ $d\nu
_{s}^{\prime }$ with $|d\mu ^{\prime }|$ $\perp $ $d\nu _{s}^{\prime }.$
Then there exists $E^{\prime }$ on which $|d\mu ^{\prime }|$ (hence $d\mu
^{\prime })$ is concentrated while $d\nu _{s}^{\prime }$ is concentrated on (%
$E^{\prime })^{c}$ $:=$ $X\backslash E^{\prime }.$ Let $E_{1}^{\prime }$ $:=$
$E^{\prime }$ $\cap $ $\{A=0\}$ and $E_{2}^{\prime }$ $:=$ $E^{\prime }$ $%
\cap $ $\{A\neq 0\}$ (note that $A$ is defined modulo a $|d\mu ^{\prime }|$%
-measure zero set in $E^{\prime }$). It follows that $d\mu $ is concentrated
on $E$ $:=$ $E_{2}^{\prime }$ $\cup $ ($E^{\prime })^{c}$. We extend the
domain of $N_{\mu }$ ($N_{\mu ^{\prime }},$ resp.) and define $N_{\mu }$ ($%
N_{\mu ^{\prime }},$ resp.) to be $0$ on $E^{c}$ $:=$ $X\backslash E$ (($%
E^{\prime })^{c},$ resp.). Let $\chi _{E}$ ($\chi _{E^{\prime }},$ resp.)
denote the characteristic function of $E$ ($E^{\prime },$ resp.), i.e., $%
\chi _{E}$ $=$ $1$ on $E$ and $\chi _{E}$ $=$ $0$ on $E^{c}.$ Following a
similar computation in (\ref{eqn3.8.2}), we then have%
\begin{eqnarray}
(N_{\mu }-N_{\mu ^{\prime }})\cdot (d\mu -d\mu ^{\prime }) &=&(1-N_{\mu
}\cdot N_{\mu ^{\prime }})(|d\mu |+|d\mu ^{\prime }|)  \label{3.8.3} \\
&=&\frac{1}{\chi _{E}+\chi _{E^{\prime }}}|N_{\mu }-N_{\mu ^{\prime
}}|^{2}(|d\mu |+|d\mu ^{\prime }|).  \notag
\end{eqnarray}

\noindent Note that since $E\cup E^{\prime }$ $=$ $X$, we have $\chi _{E}$ $%
+ $ $\chi _{E^{\prime }}$ $\neq $ $0$ on $X.$

For $\varepsilon $\ regular there holds

\begin{equation}
\mathcal{F}^{\prime }(\varepsilon )=\frac{d\mathcal{F(}d\mu _{\varepsilon })%
}{d\varepsilon }=\int_{X}N_{\varepsilon }\cdot A_{\varepsilon }|d\mu
_{\varepsilon }|.  \label{eqn3.9}
\end{equation}

\noindent Since for $\varepsilon $\textit{\ }regular we have $d\nu
_{s}^{\varepsilon }$ $=$ $0$ and hence (\ref{eqn3.9}) follows from (\ref%
{eqn3.6}). Now let $\varepsilon _{1},$ $\varepsilon _{2}$ be regular and $%
\varepsilon _{2}$ $>$ $\varepsilon _{1}$. Observe that $d\nu $ $=$ $%
A_{\varepsilon _{2}}|d\mu _{\varepsilon _{2}}|$ $=$ $A_{\varepsilon
_{1}}|d\mu _{\varepsilon _{1}}|$ and hence from (\ref{eqn3.9}) we have

\begin{eqnarray}
\mathcal{F}^{\prime }(\varepsilon _{2})-\mathcal{F}^{\prime }(\varepsilon
_{1}) &=&\int_{X}\{N_{\varepsilon _{2}}\cdot A_{\varepsilon _{2}}|d\mu
_{\varepsilon _{2}}|-N_{\varepsilon _{1}}\cdot A_{\varepsilon _{1}}|d\mu
_{\varepsilon _{1}}|\}  \label{eqn3.10} \\
&=&\int_{X}(N_{\varepsilon _{2}}-N_{\varepsilon _{1}})\cdot d\nu .  \notag
\end{eqnarray}

On the other hand, by Lemma 3.2 and Lemma 3.3 with $d\mu $ $=$ $d\mu
_{\varepsilon _{2}},$ $d\mu ^{\prime }$ $=$ $d\mu _{\varepsilon _{1}},$ (\ref%
{eqn3.8.1}) reads%
\begin{equation}
(N_{\varepsilon _{2}}-N_{\varepsilon _{1}})\cdot (d\mu _{\varepsilon
_{2}}-d\mu _{\varepsilon _{1}})=\frac{1}{2}|N_{\varepsilon
_{2}}-N_{\varepsilon _{1}}|^{2}(|d\mu _{\varepsilon _{2}}|+|d\mu
_{\varepsilon _{1}}|).  \label{eqn3.11}
\end{equation}

\noindent In view of $d\nu $ $=$ $(\varepsilon _{2}-\varepsilon _{1})^{-1}$ $%
(d\mu _{\varepsilon _{2}}-d\mu _{\varepsilon _{1}}),$ we have

\begin{equation}
\mathcal{F}^{\prime }(\varepsilon _{2})-\mathcal{F}^{\prime }(\varepsilon
_{1})=\int_{X}\frac{1}{2(\varepsilon _{2}-\varepsilon _{1})}|N_{\varepsilon
_{2}}-N_{\varepsilon _{1}}|^{2}(|d\mu _{\varepsilon _{2}}|+|d\mu
_{\varepsilon _{1}}|)\geq 0  \label{eqn3.11.0}
\end{equation}

\noindent by (\ref{eqn3.10}) and (\ref{eqn3.11}).

We remark that (\ref{eqn3.9}) generalizes Lemma 3.1 in \cite{chy}. For an
arbitrary $\varepsilon $ (regular or singular), we write $\mathcal{F}%
_{+}^{\prime }(\varepsilon )$ for $\mathcal{F}^{\prime }(\varepsilon
+)\equiv \lim_{\tilde{\varepsilon}\rightarrow \varepsilon +}\frac{\mathcal{F(%
\tilde{\varepsilon})-F(\varepsilon )}}{\tilde{\varepsilon}-\varepsilon },$
the right derivative of $\mathcal{F}$ at $\varepsilon .$ Similarly we write $%
\mathcal{F}_{-}^{\prime }(\varepsilon )$ for the left derivative $\mathcal{F}%
^{\prime }(\varepsilon -).$ Both $\mathcal{F}_{+}^{\prime }(\varepsilon )$
and $\mathcal{F}_{-}^{\prime }(\varepsilon )$ exist in view of (\ref{eqn3.6}%
). When $\varepsilon $ is regular, $\mathcal{F}_{+}^{\prime }(\varepsilon )$ 
$=$ $\mathcal{F}_{-}^{\prime }(\varepsilon )$ $=$ $\mathcal{F}^{\prime
}(\varepsilon )$ (see (\ref{eqn3.9})). We study the left and right
continuity of $\mathcal{F}_{+}^{\prime }(\varepsilon )$ and $\mathcal{F}%
_{-}^{\prime }(\varepsilon ).$

\bigskip

\textbf{Theorem 3.4.} \textit{For }$\varepsilon _{2}$\textit{\ }$>$\textit{\ 
}$\varepsilon _{1},$\textit{\ we have}%
\begin{equation}
\mathcal{F}_{+}^{\prime }(\varepsilon _{2})\geq \mathcal{F}_{-}^{\prime
}(\varepsilon _{2})\geq \mathcal{F}_{+}^{\prime }(\varepsilon _{1})\geq 
\mathcal{F}_{-}^{\prime }(\varepsilon _{1}).  \label{3.11.1}
\end{equation}%
$\mathit{\noindent }$\textit{In particular,} $\mathcal{F}^{\prime
}(\varepsilon )$\textit{\ is an increasing function of }$\varepsilon $%
\textit{\ for }$\varepsilon $\textit{\ regular. We also have the following
limits:}

\begin{eqnarray}
\lim_{\varepsilon _{2}\rightarrow \varepsilon _{1}+}\mathcal{F}_{+}^{\prime
}(\varepsilon _{2}) &=&\mathcal{F}_{+}^{\prime }(\varepsilon _{1}),\text{ }%
\lim_{\varepsilon _{2}\rightarrow \varepsilon _{1}-}\mathcal{F}_{+}^{\prime
}(\varepsilon _{2})=\mathcal{F}_{-}^{\prime }(\varepsilon _{1}),
\label{eqn3.12} \\
\lim_{\varepsilon _{2}\rightarrow \varepsilon _{1}+}\mathcal{F}_{-}^{\prime
}(\varepsilon _{2}) &=&\mathcal{F}_{+}^{\prime }(\varepsilon _{1}),\text{ }%
\lim_{\varepsilon _{2}\rightarrow \varepsilon _{1}-}\mathcal{F}_{-}^{\prime
}(\varepsilon _{2})=\mathcal{F}_{-}^{\prime }(\varepsilon _{1}).  \notag
\end{eqnarray}%
$\mathit{\noindent }$\textit{Moreover, }$\mathcal{F}$\textit{\ is convex.}

\bigskip

\textbf{Proof.} That $\mathcal{F}^{\prime }(\varepsilon )$\ is an increasing
function of $\varepsilon $\ for $\varepsilon $\ regular follows from (\ref%
{eqn3.11.0}). From (\ref{eqn3.2}) we have

\begin{eqnarray}
d\mu _{\varepsilon _{2}} &=&d\mu _{\varepsilon _{1}}+(\varepsilon
_{2}-\varepsilon _{1})d\nu  \label{eqn3.13} \\
&=&((\varepsilon _{2}-\varepsilon _{1})A_{\varepsilon _{1}}+N_{\varepsilon
_{1}})|d\mu _{\varepsilon _{1}}|+(\varepsilon _{2}-\varepsilon
_{1})N_{s}^{\varepsilon _{1}}|d\nu _{s}^{\varepsilon _{1}}|.  \notag
\end{eqnarray}%
$\mathit{\noindent }$Here we have written $d\nu _{s}^{\varepsilon _{1}}$ $=$ 
$N_{s}^{\varepsilon _{1}}|d\nu _{s}^{\varepsilon _{1}}|.$ We need the
following lemma to describe $N_{\varepsilon _{2}}.$

\bigskip

\textbf{Lemma 3.5.} \textit{Let }$d\lambda ,$\textit{\ }$d\tau ,$\textit{\
and }$d\rho $\textit{\ be bounded vector-valued measures.} \textit{Suppose }$%
d\lambda $\textit{\ }$=$\textit{\ }$\vec{B}|d\tau |$\textit{\ }$+$\textit{\ }%
$\vec{C}|d\rho |,$\textit{\ }$|d\tau |$\textit{\ }$\perp $\textit{\ }$|d\rho
|,$\textit{\ and }$|d\lambda |$\textit{\ }$\ll $\textit{\ }$|d\tau |$\textit{%
\ }$+$\textit{\ }$|d\rho |$\textit{\ }$\ll $\textit{\ }$|d\lambda |.$\textit{%
\ Then }$\vec{B}$\textit{\ }$\neq $\textit{\ }$0$\textit{\ a.e. }$[|d\tau
|], $\textit{\ }$\vec{C}$\textit{\ }$\neq $\textit{\ }$0$\textit{\ a.e. }$%
[|d\rho |],$\textit{\ and }%
\begin{equation}
\vec{N}_{\lambda }=\frac{\vec{B}}{|\vec{B}|}\text{ a.e. }[|d\tau |];\text{ }=%
\frac{\vec{C}}{|\vec{C}|}\text{ a.e. }[|d\rho |]  \label{eqn3.14}
\end{equation}

\noindent \textit{where we write }$d\lambda $\textit{\ }$=$\textit{\ }$\vec{N%
}_{\lambda }|d\lambda |.$

\bigskip

%TCIMACRO{\TeXButton{Proof}{\proof} }%
%BeginExpansion
\proof
%EndExpansion
Since $|d\tau |$\textit{\ }$\ll $\textit{\ }$|d\lambda |$ ($|d\rho |$\textit{%
\ }$\ll $\textit{\ }$|d\lambda |,$ respectively), we can find a
vector-valued function $\vec{h}_{\tau }$ $\in $ $L^{1}(|d\lambda |)$ ($\vec{h%
}_{\rho }$ $\in $ $L^{1}(|d\lambda |)$, respectively) such that 
\begin{equation*}
d\tau =\vec{h}_{\tau }|d\lambda |\text{ (}d\rho =\vec{h}_{\rho }|d\lambda |,%
\text{ respectively).}
\end{equation*}

\noindent It follows that $|d\tau |$ $=$ $|\vec{h}_{\tau }||d\lambda |$ $=$ $%
|\vec{h}_{\tau }|(|\vec{B}||d\tau |$\textit{\ }$+$\textit{\ }$|\vec{C}%
||d\rho |).$ So $|\vec{h}_{\tau }||\vec{B}|$ $=$ $1$ a.e. $[|d\tau |]$ and $|%
\vec{h}_{\tau }||\vec{C}|$ $=$ $0$ a.e. $[|d\rho |].$ Therefore $\vec{B}$%
\textit{\ }$\neq $\textit{\ }$0$\textit{\ }a.e\textit{. }$[|d\tau |].$
Similarly we have $\vec{C}$\textit{\ }$\neq $\textit{\ }$0$\textit{\ }a.e.%
\textit{\ }$[|d\rho |]$ and hence

\begin{eqnarray}
\vec{h}_{\tau } &=&0\text{ a.e. }[|d\rho |],\text{ }|\vec{h}_{\tau }|=\frac{1%
}{|\vec{B}|}\text{ a.e. }[|d\tau |],\text{ and also}  \label{eqn3.15} \\
\vec{h}_{\rho } &=&0\text{ a.e. }[|d\tau |],\text{ }|\vec{h}_{\rho }|=\frac{1%
}{|\vec{C}|}\text{ a.e. }[|d\rho |]  \notag
\end{eqnarray}%
\noindent by symmetry. Now we compute $\vec{N}_{\lambda }|d\lambda |$ $=$ $%
d\lambda $ $=$ $\vec{B}|d\tau |$\textit{\ }$+$\textit{\ }$\vec{C}|d\rho |$ $%
= $ $\vec{B}|\vec{h}_{\tau }||d\lambda |$ $+$ $\vec{C}|\vec{h}_{\rho
}||d\lambda |$ $=$ $(\vec{B}|\vec{h}_{\tau }|$ $+$ $\vec{C}|\vec{h}_{\rho
}|)|d\lambda |.$ It then follows that%
\begin{equation}
\vec{N}_{\lambda }=\vec{B}|\vec{h}_{\tau }|+\vec{C}|\vec{h}_{\rho }|\text{
a.e. }[|d\lambda |].  \label{eqn3.16}
\end{equation}

\noindent Since $|d\tau |$\textit{\ }$\perp $\textit{\ }$|d\rho |$, we
obtain (\ref{eqn3.14}) from (\ref{eqn3.16}) in view of (\ref{eqn3.15}).

%TCIMACRO{\TeXButton{End Proof}{\endproof}}%
%BeginExpansion
\endproof%
%EndExpansion

\bigskip

%TCIMACRO{\TeXButton{Proof}{\proof} }%
%BeginExpansion
\proof
%EndExpansion
\textbf{(of Theorem 3.4 continued)}

From Lemma 3.5 we express $N_{\varepsilon _{2}}$ as follows:

\begin{equation}
N_{\varepsilon _{2}}=\left\{ 
\begin{array}{l}
\frac{(\varepsilon _{2}-\varepsilon _{1})A_{\varepsilon _{1}}+N_{\varepsilon
_{1}}}{|(\varepsilon _{2}-\varepsilon _{1})A_{\varepsilon
_{1}}+N_{\varepsilon _{1}}|}\text{ \ a.e. }[|d\mu _{\varepsilon _{1}}|] \\ 
\frac{\varepsilon _{2}-\varepsilon _{1}}{|\varepsilon _{2}-\varepsilon _{1}|}%
N_{s}^{\varepsilon _{1}}\text{ \ a.e. }[|d\nu _{s}^{\varepsilon _{1}}|].%
\end{array}%
\right.  \label{eqn3.17}
\end{equation}

\noindent Now for $\varepsilon _{2}$ regular ($\varepsilon _{1}$ may not be
regular) we compute%
\begin{eqnarray}
\mathcal{F}^{\prime }(\varepsilon _{2}) &=&\int_{X}N_{\varepsilon _{2}}\cdot
d\nu  \label{eqn3.18} \\
&=&\int_{X}N_{\varepsilon _{2}}\cdot (A_{\varepsilon _{1}}|d\mu
_{\varepsilon _{1}}|+d\nu _{s}^{\varepsilon _{1}})\text{ \ (by (\ref{eqn3.2}%
))}  \notag \\
&=&\int_{X}\frac{(\varepsilon _{2}-\varepsilon _{1})A_{\varepsilon
_{1}}+N_{\varepsilon _{1}}}{|(\varepsilon _{2}-\varepsilon
_{1})A_{\varepsilon _{1}}+N_{\varepsilon _{1}}|}\cdot A_{\varepsilon
_{1}}|d\mu _{\varepsilon _{1}}|+\frac{\varepsilon _{2}-\varepsilon _{1}}{%
|\varepsilon _{2}-\varepsilon _{1}|}N_{s}^{\varepsilon _{1}}\cdot d\nu
_{s}^{\varepsilon _{1}}  \notag
\end{eqnarray}

\noindent by (\ref{eqn3.17}). Observe that $N_{s}^{\varepsilon _{1}}\cdot
d\nu _{s}^{\varepsilon _{1}}$ $=$ $|d\nu _{s}^{\varepsilon _{1}}|$ and the
integrand in (\ref{eqn3.18}) is bounded by $|A_{\varepsilon _{1}}||d\mu
_{\varepsilon _{1}}|$ $+$ $|d\nu _{s}^{\varepsilon _{1}}|$ (which is
independent of $\varepsilon _{2}$ and integrable by assumption). We can
therefore apply the Lebesgue dominated convergence theorem to get

\begin{eqnarray}
\lim_{\varepsilon _{2}\rightarrow \varepsilon _{1}\pm }\mathcal{F}^{\prime
}(\varepsilon _{2}) &=&\int_{X}N_{\varepsilon _{1}}\cdot A_{\varepsilon
_{1}}|d\mu _{\varepsilon _{1}}|\pm |d\nu _{s}^{\varepsilon _{1}}|
\label{eqn3.19} \\
&=&\mathcal{F}_{\pm }^{\prime }(\varepsilon _{1})  \notag
\end{eqnarray}

\noindent by (\ref{eqn3.6}). Since $\mathcal{F}^{\prime }(\varepsilon )$ is
increasing for $\varepsilon $ regular and the set of regular values is
dense, we can easily deduce (\ref{3.11.1}) from (\ref{eqn3.19}). Thus we
have $\lim_{\varepsilon _{2}\rightarrow \varepsilon _{1}+}\mathcal{F}%
_{+}^{\prime }(\varepsilon _{2})$ $=$ $\mathcal{F}_{+}^{\prime }(\varepsilon
_{1})$ and $\lim_{\varepsilon _{2}\rightarrow \varepsilon _{1}+}\mathcal{F}%
_{-}^{\prime }(\varepsilon _{2})$ $=$ $\mathcal{F}_{+}^{\prime }(\varepsilon
_{1}).$ Similarly we also have $\lim_{\varepsilon _{2}\rightarrow
\varepsilon _{1}-}\mathcal{F}_{+}^{\prime }(\varepsilon _{2})$ $=$ $\mathcal{%
F}_{-}^{\prime }(\varepsilon _{1})$ and $\lim_{\varepsilon _{2}\rightarrow
\varepsilon _{1}-}\mathcal{F}_{-}^{\prime }(\varepsilon _{2})$ $=$ $\mathcal{%
F}_{-}^{\prime }(\varepsilon _{1}).$ We have proved (\ref{eqn3.12}). That $%
\mathcal{F}$ is convex follows from (\ref{3.11.1}) by elementary calculus.

%TCIMACRO{\TeXButton{End Proof}{\endproof}}%
%BeginExpansion
\endproof%
%EndExpansion

\bigskip

We remark that Theorem 3.4 generalizes Lemma 3.2 in \cite{chy}.

\bigskip

\textbf{Example 3.1. }Consider a $C^{2}$ smooth graph $\Sigma $ $=$ $%
\{(x_{1},$ $x_{2},$ $...,$ $x_{m},$ $u(x_{1},$ $x_{2},$ $...,$ $x_{m}))\}$
in $R^{m+1}.$ Let $d\mu $ $:=$ $(u_{x_{1}},$ $u_{x_{2}},$ $...,$ $u_{x_{m}},$
$-1)$ $d^{m}x$ where we recall that $d^{m}x$ :$=$ $dx_{1}$ $\wedge $ $dx_{2}$
$\wedge $ $...$ $\wedge $ $dx_{m}$ be the $R^{m+1}$-valued measure defined
on a bounded domain $X$ $\subset $ $R^{m},$ associated with the Euclidean
normal to $\Sigma .$ Then $|d\mu |$ $=$ $\sqrt{%
1+u_{x_{1}}^{2}+...+u_{x_{m}}^{2}}$ $d^{m}x$ is the area element of $\Sigma $
with respect to the metric induced from the Euclidean metric on $R^{m+1}.$
Let $d\nu $ $\equiv $ $(v_{x_{1}},$ $v_{x_{2}},$ $...,$ $v_{x_{m}},$ $0)$ $%
d^{m}x$ where $v$ $\in $ $C_{0}^{\infty }(X).$ So from $d\mu _{\varepsilon }$
$=$ $N_{\varepsilon }$ $|d\mu _{\varepsilon }|$, $d\nu $ $=$ $A_{\varepsilon
}$ $|d\mu _{\varepsilon }|,$ and $|d\mu _{\varepsilon }|$ $=$ $\sqrt{%
1+|\nabla u+\varepsilon \nabla v|^{2}}$ $d^{m}x$ ($d\nu _{s}^{\varepsilon }$ 
$=$ $0$ since $|d\mu _{\varepsilon }|$ is strictly positive; so each $%
\varepsilon $ is regular), we obtain

\begin{equation}
N_{\varepsilon }=\frac{(\nabla u+\varepsilon \nabla v,-1)}{\sqrt{1+|\nabla
u+\varepsilon \nabla v|^{2}}},A_{\varepsilon }=\frac{(\nabla v,0)}{\sqrt{%
1+|\nabla u+\varepsilon \nabla v|^{2}}}  \label{eqn3.21}
\end{equation}

\noindent where $\nabla $ denotes the gradient in $R^{m}.$ By Theorem B we
have the first variation of the area $\mathcal{F(}0)$ $=$ $\mathcal{F(}d\mu
) $ of $\Sigma :$

\begin{eqnarray}
\mathcal{F}^{\prime }\mathcal{(}0) &=&\int_{X}N_{0}\cdot A_{0}|d\mu |
\label{eqn3.22} \\
&=&\int_{X}\frac{\nabla u\cdot \nabla v}{\sqrt{1+|\nabla u|^{2}}}d^{m}x 
\notag \\
&=&-\int_{X}div(\frac{\nabla u}{\sqrt{1+|\nabla u|^{2}}})vd^{m}x  \notag
\end{eqnarray}

\noindent by (\ref{eqn3.21}) and the divergence theorem. Notice that $div(%
\frac{\nabla u}{\sqrt{1+|\nabla u|^{2}}})$ in (\ref{eqn3.22}) is the
(Riemannian) mean curvature of $\Sigma $ in $R^{m+1}.$ We have recovered the
classical first variation formula for the area of a graph in the Euclidean
space.

\bigskip

\textbf{Example 3.2. }Consider a $C^{1}$ smooth graph $\Sigma $ $=$ $%
\{(x_{1},$ $x_{1^{\prime }},$ $...,$ $x_{n},$ $x_{n^{\prime }},$ $u(x_{1},$ $%
x_{1^{\prime }},$ $...,$ $x_{n},$ $x_{n^{\prime }}))\}$ in the Heisenberg
group viewed as $R^{2n+1}$ with the standard flat pseudohermitian structure
(see \cite{chmy}). Recall that $\vec{X}^{\ast }$ $=$ $(x_{1^{\prime }},$ $%
-x_{1},$ $x_{2^{\prime }},$ $-x_{2},...,$ $x_{n^{\prime }},$ $-x_{n}).$ Let $%
\nabla $ denote the gradient opeator in $R^{2n}.$ Let $d\mu $ $:=$ $(\nabla
u $ $-$ $\vec{X}^{\ast })$ $d^{2n}x$ where $d^{2n}x$ $:=$ $dx_{1}$ $\wedge $ 
$dx_{1^{\prime }}$ $\wedge $ $...$ $\wedge $ $dx_{n}$ $\wedge $ $%
dx_{n^{\prime }}$ and $d\nu $ $\equiv $ $(\nabla \varphi )$ $d^{2n}x$ be two 
$R^{2n}$-valued measures defined on a bounded domain $\Omega $ $\subset $ $%
R^{2n}$ $(\varphi $ $\in $ $C_{0}^{1}(\Omega ),$ say$).$ So $|d\mu |$ $=$ $%
|\nabla u$ $-$ $\vec{X}^{\ast }|d^{2n}x$ is the $p$-area element. Denote the
singular set $\{\nabla u$ $-$ $\vec{X}^{\ast }$ $=$ $0\}$ by $S(u).$ Write $%
d\mu $ $=$ $N_{0}$ $|d\mu |\ $and $d\nu $ $=$ $A_{0}$ $|d\mu |$ $+$ $d\nu
_{s}$ where

\begin{eqnarray}
N_{0} &=&\frac{\nabla u-\vec{X}^{\ast }}{|\nabla u-\vec{X}^{\ast }|},\text{ }%
A_{0}=\frac{\nabla \varphi }{|\nabla u-\vec{X}^{\ast }|}\text{ on }\Omega
\backslash S(u)  \label{eqn3.23} \\
d\nu _{s} &=&(\nabla \varphi )d^{2n}x\text{ on }S(u).  \notag
\end{eqnarray}

\noindent Note that $|d\mu |$ is concentrated on $\Omega \backslash S(u)$
while $d\nu _{s}$ is concentrated on $S(u).$ By Theorem B and (\ref{eqn3.23}%
) we have the first variation of the $p$-area $\mathcal{F(}0)$ $=$ $\mathcal{%
F(}d\mu )$ of $\Sigma :$%
\begin{eqnarray}
\mathcal{F}^{\prime }\mathcal{(}0\pm ) &=&\int_{\Omega }N_{0}\cdot
A_{0}|d\mu |\text{ }\pm \text{ }|d\nu _{s}|  \label{eqn3.24} \\
&=&\int_{\Omega \backslash S(u)}\frac{(\nabla u-\vec{X}^{\ast })\cdot \nabla
\varphi }{|\nabla u-\vec{X}^{\ast }|}d^{2n}x\pm \int_{S(u)}|\nabla \varphi
|d^{2n}x  \notag \\
&=&\int_{\Omega \backslash S(u)}\func{div}(\varphi \frac{\nabla u-\vec{X}%
^{\ast }}{|\nabla u-\vec{X}^{\ast }|})d^{2n}x  \notag \\
&&-\int_{\Omega \backslash S(u)}\varphi \func{div}(\frac{\nabla u-\vec{X}%
^{\ast }}{|\nabla u-\vec{X}^{\ast }|})d^{2n}x\pm \int_{S(u)}|\nabla \varphi
|d^{2n}x  \notag
\end{eqnarray}

\noindent (cf. (3.3) in \cite{chy}). We remark that the Lebesgue measure of $%
S(u)$ vanishes for $u$ $\in $ $C^{2}$ (in this case, compare (\ref{eqn3.24})
with the first variation formula in \cite{Ri2}) or $C^{1,1}$ while there
exists $u$ $\in $ $\cap _{0<\alpha <1}C^{1,\alpha }$ such that $S(u)$ has
positive Lebesgue measure according to Balogh (\cite{Ba}).

\bigskip

\textbf{Example 3.3.} For basic material in this example, the readers are
referred to (\cite{chmy}). Let $(M,$ $J,$ $\Theta )$ be a 3-dimensional
oriented pseudohermitian manifold. Consider a $C^{2}$ smooth orientable
surface $\Sigma $ $\subset $ $M.$ Let $\xi $ $\equiv $ $\ker \Theta $ denote
the contact bundle. Let $e_{1}$ $\in $ $T\Sigma $ $\cap $ $\xi $ denote a
characteristic vector of unit length with respect to the Levi metric $G$ $=$ 
$\frac{1}{2}d\Theta (\cdot ,J\cdot )$ (at a nonsingular point)$.$ Let $e_{2}$
$\equiv $ $Je_{1}$ and $T$ denote the Reeb vector field associated to $%
\Theta .$ Let $\{e^{1},$ $e^{2},$ $\Theta \}$ denote the coframe field dual
to the frame field $\{e_{1},$ $e_{2},$ $T\}.$ The adapted (or left
invariant) metric on $M$ is defined by $h$ $=$ $\Theta ^{2}$ $+$ $G$ $=$ $%
\Theta ^{2}$ $+$ $(e^{1})^{2}$ $+$ $(e^{2})^{2}$ (if restricted on the
nonsingular domain)$.$ It follows that

\begin{equation}
\tilde{e}_{1}=e_{1},\tilde{e}_{2}=-\frac{\alpha e_{2}+T}{\sqrt{1+\alpha ^{2}}%
},N=\frac{e_{2}-\alpha T}{\sqrt{1+\alpha ^{2}}}  \label{eqn3.25}
\end{equation}

\noindent form an orthonormal basis with respect to $h$ (recall that $\alpha 
$ is defined so that $\alpha e_{2}+T$ $\in $ $T\Sigma ).$ Denote the
projection of the unit normal $N$ onto $\xi $ by $N_{\xi }.$ Denote the
Riemannian area element of $\Sigma $ induced from $h$ by $d\Sigma .$ Let $%
\tilde{e}^{1}$ $=$ $e^{1},$ $\tilde{e}^{2}$ $=$ $-\frac{\alpha e^{2}+\Theta 
}{\sqrt{1+\alpha ^{2}}},$ and $\tilde{e}^{3}$ $=$ $\frac{e^{2}-\alpha \Theta 
}{\sqrt{1+\alpha ^{2}}}$ be the coframe field dual to $\tilde{e}_{1},$ $%
\tilde{e}_{2},$ $N$ in (\ref{eqn3.25}). We have

\begin{equation}
N_{\xi }=\frac{e_{2}}{\sqrt{1+\alpha ^{2}}},d\Sigma =\tilde{e}^{1}\wedge 
\tilde{e}^{2}  \label{eqn3.26}
\end{equation}

\noindent (assuming that $\Sigma $ is oriented so that the second equality
in (\ref{eqn3.26}) holds). Let $|$ $\cdot $ $|_{h}$ denote the length with
respect to the metric $h.$ From (\ref{eqn3.26}) we can now compute

\begin{eqnarray}
|N_{\xi }|_{h}d\Sigma &=&\frac{1}{\sqrt{1+\alpha ^{2}}}\tilde{e}^{1}\wedge 
\tilde{e}^{2}  \label{eqn3.27} \\
&=&\frac{-1}{\sqrt{1+\alpha ^{2}}}e^{1}\wedge \frac{\alpha e^{2}+\Theta }{%
\sqrt{1+\alpha ^{2}}}  \notag \\
&=&\Theta \wedge e^{1}  \notag
\end{eqnarray}

\noindent (on the nonsingular domain; $=$ $0$ on the singular set) by noting
that $e^{1}$ $\wedge $ $e^{2}$ $=$ $\alpha e^{1}$ $\wedge $ $\Theta $ on $%
\Sigma .$ For $M$ being the Heisenberg group, (\ref{eqn3.27}) was pointed
out in \cite{Ri2}. So we learn from (\ref{eqn3.27}) that the general $p$%
-area element can also be viewed as the total variation measure of a $TM$ or 
$\xi $-valued measure $N_{\xi }d\Sigma $ on $\Sigma .$

By the way we will compute the first variation formula for variations having
support containing the singular set (in \cite{chmy} we computed it for
variations having support away from the singular set). For simplicity we
assume that $\Sigma $ is $C^{1}$ smooth, oriented, and $\Sigma \backslash
S_{\Sigma }$ is $C^{2},$ where $S_{\Sigma }$ denotes the singular set
consisting of a $C^{1}$ smooth curve. Suppose $S_{\Sigma }$ divides $\Sigma $
into two pieces with boundaries $S_{\Sigma }^{+},$ $S_{\Sigma }^{-}$
reversely oriented on $S_{\Sigma }.$ Let $v$ be a $C^{\infty }$ smooth
vector field of $M$ with support away from $\partial \Sigma $ when
restricted to $\Sigma .$ We write $v$ $=$ $v_{1}e_{1}$ $+$ $v_{2}e_{2}$ $+$ $%
fT$ (in nonsingular region)$.$ Compute the variation of the general $p$-area
in the direction $v:$

\begin{eqnarray}
\delta _{v}\int_{\Sigma }\Theta \wedge e^{1} &=&\int_{\Sigma \backslash
S_{\Sigma }}L_{v}(\Theta \wedge e^{1})  \label{eqn3.28} \\
&=&\int_{\Sigma \backslash S_{\Sigma }}d\circ i_{v}(\Theta \wedge
e^{1})+i_{v}\circ d(\Theta \wedge e^{1})  \notag \\
&=&(\int_{S_{\Sigma }^{+}}+\int_{S_{\Sigma }^{-}})(fe^{1}-v_{1}\Theta
)+\int_{\Sigma \backslash S_{\Sigma }}(f\alpha -v_{2})H\Theta \wedge e^{1} 
\notag
\end{eqnarray}

\noindent by (2.8') in \cite{chmy}, where $H$ denotes the $p$-mean curvature
of $\Sigma $. We say that $\Sigma $ is stationary if $\delta
_{v}\int_{\Sigma }\Theta \wedge e^{1}$ $=$ $0$ for all $v.$ Then by (\ref%
{eqn3.28}) and $\Theta $ $=$ $0$ on $S_{\Sigma },$ we learn that if $\Sigma $
is stationary, then $H$ $=$ $0$ by taking $v$ with support away from $%
S_{\Sigma },$ and hence there holds

\begin{equation}
\int_{S_{\Sigma }^{+}}fe^{1}+\int_{S_{\Sigma }^{-}}fe^{1}=0.  \label{eqn3.29}
\end{equation}

\noindent Let $\tau $ denote the positive unit vector tangent to $S_{\Sigma
}^{+}.$ Assume that we can extend $e_{1}$ continuously to $S_{\Sigma }$ from
both sides. Denote the extensions of $e_{1}$ and $e^{1}$ on $S_{\Sigma }^{+}$
($S_{\Sigma }^{-},$ respectively) by $e_{1}^{+}$ and $e_{+}^{1}$ ($e_{1}^{-}$
and $e_{-}^{1},$ respectively). Then from (\ref{eqn3.29}) we have

\begin{equation}
0=e_{+}^{1}(\tau )+e_{-}^{1}(-\tau )=e_{1}^{+}\cdot \tau -e_{1}^{-}\cdot \tau
\label{eqn3.30}
\end{equation}

\noindent where \textquotedblright $\cdot "$ denote the inner product with
respect to the adapted metric $h$ or the Levi metric $G$ (note that $%
e_{1}^{+},$ $e_{1}^{-},$ and $\tau $ are all in $\xi ).$ (\ref{eqn3.30}) is
the "incident angle = reflected angle" condition on the singular curves for
a $p$-area stationary surface. When $\Sigma $ is $C^{2}$ (including the
singular set $S_{\Sigma }),$ both "angles" must be 90 degrees, i.e., $%
e_{1}^{+}\cdot \tau $ $=$ $e_{1}^{-}\cdot \tau $ $=$ $0$ since $e_{1}^{+}$ $%
= $ $-e_{1}^{-}$ according to (generalized) Proposition 3.5 in \cite{chmy}
(see a remark in Section 7 for generalizing the results in Section 3). We
studied condition (\ref{eqn3.30}) for a (generalized) stationary graph in
the Heisenberg group (see Theorem 6.3 in \cite{chy}). Ritor\'{e} and Rosales
(\cite{Ri2}) obtained the same result for a $C^{2}$ smooth, oriented
(immersed) surface in the Heisenberg group.

\bigskip

\section{\textbf{Second variation and proof of Theorem C}}

Recall that for $\varepsilon $ regular we have (cf. (\ref{eqn3.9}))

\begin{equation}
\mathcal{F}^{\prime }\mathcal{(}\varepsilon )=\frac{d\mathcal{F(}d\mu
_{\varepsilon })}{d\varepsilon }\mathcal{=}\int_{X}N_{\varepsilon }\cdot
d\nu .  \label{eqn4.1}
\end{equation}

\noindent First from (\ref{eqn4.1}), we want to compute 
\begin{equation}
\lim_{\varepsilon _{2}\rightarrow \varepsilon _{1}}\frac{\mathcal{F}^{\prime
}\mathcal{(}\varepsilon _{2})-\mathcal{F}^{\prime }\mathcal{(}\varepsilon
_{1})}{\varepsilon _{2}-\varepsilon _{1}}=\lim_{\varepsilon _{2}\rightarrow
\varepsilon _{1}}\int_{X}(\frac{N_{\varepsilon _{2}}-N_{\varepsilon _{1}}}{%
\varepsilon _{2}-\varepsilon _{1}})\cdot d\nu  \label{eqn4.2}
\end{equation}

\noindent for $\varepsilon _{2}$, $\varepsilon _{1}$ regular. From (\ref%
{eqn3.13}) we have%
\begin{equation}
d\mu _{\varepsilon _{2}}=((\varepsilon _{2}-\varepsilon _{1})A_{\varepsilon
_{1}}+N_{\varepsilon _{1}})|d\mu _{\varepsilon _{1}}|  \label{eqn3.7}
\end{equation}

\noindent for $\varepsilon _{1}$ regular. Taking the absolute value (total
variation) of both sides in (\ref{eqn3.7}) gives%
\begin{equation}
|d\mu _{\varepsilon _{2}}|=|(\varepsilon _{2}-\varepsilon
_{1})A_{\varepsilon _{1}}+N_{\varepsilon _{1}}||d\mu _{\varepsilon _{1}}|
\label{eqn3.8}
\end{equation}

\noindent or

\begin{equation}
|d\mu _{\varepsilon _{1}}|=\frac{1}{|(\varepsilon _{2}-\varepsilon
_{1})A_{\varepsilon _{1}}+N_{\varepsilon _{1}}|}|d\mu _{\varepsilon _{2}}|.
\label{eqn4.5}
\end{equation}%
\noindent Note that $(\varepsilon _{2}-\varepsilon _{1})A_{\varepsilon
_{1}}+N_{\varepsilon _{1}}$ $\neq $ $0$ a.e. $[|d\mu _{\varepsilon _{1}}|]$
and $[|d\mu _{\varepsilon _{2}}|]$ as shown below.

\bigskip

\textbf{Lemma 4.1}.\textit{\ Let }$\varepsilon _{2}$\textit{\ and }$%
\varepsilon _{1}$\textit{\ be regular. Then we have }$(\varepsilon
_{2}-\varepsilon _{1})A_{\varepsilon _{1}}+N_{\varepsilon _{1}}$\textit{\ }$%
\neq $\textit{\ }$0$\textit{\ a.e. }$[|d\mu _{\varepsilon _{1}}|]$\textit{\
(and hence also a.e. }$[|d\mu _{\varepsilon _{2}}|]$ \textit{by Lemma 3.2 (2)%
}$).$

\bigskip

%TCIMACRO{\TeXButton{Proof}{\proof} }%
%BeginExpansion
\proof
%EndExpansion
Suppose there is a $|d\mu _{\varepsilon _{1}}|$-measurable set $S$ such that 
$|d\mu _{\varepsilon _{1}}|(S)$ $>$ $0$ while $(\varepsilon _{2}-\varepsilon
_{1})A_{\varepsilon _{1}}+N_{\varepsilon _{1}}$\textit{\ }$=$\textit{\ }$0.$
By (\ref{eqn3.8}) we have $|d\mu _{\varepsilon _{2}}|(S)$ $=$ $0$,
contradicting $|d\mu _{\varepsilon _{1}}|$\textit{\ }$\ll $\textit{\ }$|d\mu
_{\varepsilon _{2}}|$ as asserted in Lemma 3.2 (2).

%TCIMACRO{\TeXButton{End Proof}{\endproof}}%
%BeginExpansion
\endproof%
%EndExpansion

\bigskip

Substituting the first equality of (\ref{eqn3.2}) with $\varepsilon $ $=$ $%
\varepsilon _{2}$ and (\ref{eqn4.5}) into (\ref{eqn3.7}), we get

\begin{equation}
N_{\varepsilon _{2}}=\frac{(\varepsilon _{2}-\varepsilon _{1})A_{\varepsilon
_{1}}+N_{\varepsilon _{1}}}{|(\varepsilon _{2}-\varepsilon
_{1})A_{\varepsilon _{1}}+N_{\varepsilon _{1}}|}.  \label{eqn4.6}
\end{equation}%
\noindent From (\ref{eqn4.6}) we can write%
\begin{equation}
N_{\varepsilon _{2}}-N_{\varepsilon _{1}}=(I)+(II)  \label{eqn4.7}
\end{equation}

\noindent where 
\begin{eqnarray*}
(I) &=&\frac{(\varepsilon _{2}-\varepsilon _{1})A_{\varepsilon
_{1}}+N_{\varepsilon _{1}}}{|(\varepsilon _{2}-\varepsilon
_{1})A_{\varepsilon _{1}}+N_{\varepsilon _{1}}|}-[(\varepsilon
_{2}-\varepsilon _{1})A_{\varepsilon _{1}}+N_{\varepsilon _{1}}], \\
(II) &=&[(\varepsilon _{2}-\varepsilon _{1})A_{\varepsilon
_{1}}+N_{\varepsilon _{1}}]-N_{\varepsilon _{1}}.
\end{eqnarray*}

\noindent So we can estimate%
\begin{eqnarray}
|(II)| &=&|\varepsilon _{2}-\varepsilon _{1}||A_{\varepsilon _{1}}|\text{ \
and}  \label{eqn4.8} \\
|(I)| &=&|\frac{(\varepsilon _{2}-\varepsilon _{1})A_{\varepsilon
_{1}}+N_{\varepsilon _{1}}}{|(\varepsilon _{2}-\varepsilon
_{1})A_{\varepsilon _{1}}+N_{\varepsilon _{1}}|}(1-|(\varepsilon
_{2}-\varepsilon _{1})A_{\varepsilon _{1}}+N_{\varepsilon _{1}}|)|  \notag \\
&\leq &|\varepsilon _{2}-\varepsilon _{1}||A_{\varepsilon _{1}}|  \notag
\end{eqnarray}

\noindent by noting that $1$ $=$ $|N_{\varepsilon _{1}}|$ and making use of
the triangle inequality (a.e. for $|d\mu _{\varepsilon _{1}}|$ and also for $%
|d\mu _{\varepsilon _{2}}|$ by Lemma 4.1 (1)). From (\ref{eqn4.7}) and (\ref%
{eqn4.8}) we have

\begin{equation}
|\frac{N_{\varepsilon _{2}}-N_{\varepsilon _{1}}}{\varepsilon
_{2}-\varepsilon _{1}}|\leq 2|A_{\varepsilon _{1}}|.  \label{eqn4.9}
\end{equation}

\noindent Since $|A_{\varepsilon _{1}}||d\nu |$ $=$ $|A_{\varepsilon
_{1}}|^{2}|d\mu _{\varepsilon _{1}}|$ is integrable by assumption, we can
therefore apply the Lebesgue dominated convergence theorem to get

\begin{equation}
\lim_{\varepsilon _{2}\rightarrow \varepsilon _{1}}\int_{X}(\frac{%
N_{\varepsilon _{2}}-N_{\varepsilon _{1}}}{\varepsilon _{2}-\varepsilon _{1}}%
)\cdot d\nu =\int_{X}(\lim_{\varepsilon _{2}\rightarrow \varepsilon _{1}}%
\frac{N_{\varepsilon _{2}}-N_{\varepsilon _{1}}}{\varepsilon
_{2}-\varepsilon _{1}})\cdot d\nu  \label{eqn4.10}
\end{equation}

\noindent by (\ref{eqn4.9}). Let%
\begin{equation*}
f(t)\equiv \frac{tA_{\varepsilon _{1}}+N_{\varepsilon _{1}}}{%
|tA_{\varepsilon _{1}}+N_{\varepsilon _{1}}|}.
\end{equation*}

\noindent Recall that $tA_{\varepsilon _{1}}+N_{\varepsilon _{1}}$ $\neq $ $%
0 $ a.e. (for $|d\mu _{\varepsilon _{1}}|)$ for $t$ $=$ $\varepsilon
_{2}-\varepsilon _{1}$ and $0.$ A straightforward computation shows that

\begin{equation}
f^{\prime }(t)=\frac{A_{\varepsilon _{1}}-(A_{\varepsilon _{1}}\cdot
N_{\varepsilon _{1}})N_{\varepsilon _{1}}+t[(A_{\varepsilon _{1}}\cdot
N_{\varepsilon _{1}})A_{\varepsilon _{1}}-|A_{\varepsilon
_{1}}|^{2}N_{\varepsilon _{1}}]}{|tA_{\varepsilon _{1}}+N_{\varepsilon
_{1}}|^{3}}.  \label{eqn4.11}
\end{equation}

\noindent It follows that 
\begin{equation}
f^{\prime }(t)\cdot A_{\varepsilon _{1}}=\frac{|A_{\varepsilon
_{1}}|^{2}-|(A_{\varepsilon _{1}}\cdot N_{\varepsilon _{1}})|^{2}}{%
|tA_{\varepsilon _{1}}+N_{\varepsilon _{1}}|^{3}}\geq 0  \label{eqn4.12}
\end{equation}

\noindent by Cauchy's inequality (noting that $|N_{\varepsilon _{1}}|$ $=$ $%
1 $), and

\begin{eqnarray}
\lim_{\varepsilon _{2}\rightarrow \varepsilon _{1}}\frac{\mathcal{F}^{\prime
}\mathcal{(}\varepsilon _{2})-\mathcal{F}^{\prime }\mathcal{(}\varepsilon
_{1})}{\varepsilon _{2}-\varepsilon _{1}} &=&\int_{X}(\lim_{\varepsilon
_{2}\rightarrow \varepsilon _{1}}\frac{f(\varepsilon _{2}-\varepsilon
_{1})-f(0)}{\varepsilon _{2}-\varepsilon _{1}})\cdot d\nu  \label{eqn4.13} \\
&=&\int_{X}f^{\prime }(0)\cdot A_{\varepsilon _{1}}|d\mu _{\varepsilon _{1}}|
\notag \\
&=&\int_{X}\{|A_{\varepsilon _{1}}|^{2}-|(A_{\varepsilon _{1}}\cdot
N_{\varepsilon _{1}})|^{2}\}|d\mu _{\varepsilon _{1}}|\geq 0  \notag
\end{eqnarray}

\noindent by (\ref{eqn4.2}), (\ref{eqn4.10}), (\ref{eqn3.2}) (with $d\nu
_{s}^{\varepsilon _{1}}$ $=$ $0),$ and (\ref{eqn4.12}) for $\varepsilon
_{2}, $ $\varepsilon _{1}$ regular. We have proved Theorem C (1) (\ref%
{eqn1.6}).

Next we are going to prove Theorem C (2). Take arbitrary $\varepsilon _{1},$ 
$\varepsilon _{2},$ $\varepsilon _{1}$ $\neq $ $\varepsilon _{2}.$ First we
want to express $\mathcal{F}_{\pm }^{\prime }\mathcal{(}\varepsilon _{2})$
in terms of $|d\mu _{\varepsilon _{1}}|$ and $|d\nu _{s}^{\varepsilon
_{1}}|. $ Since $|d\mu _{\varepsilon _{1}}|$ $\perp $ $|d\nu
_{s}^{\varepsilon _{1}}|,$ there exists $E_{\varepsilon _{1}}$ such that $%
|d\mu _{\varepsilon _{1}}|$ is concentrated on $E_{\varepsilon _{1}}$ while $%
|d\nu _{s}^{\varepsilon _{1}}|$ is concentrated on $E_{\varepsilon _{1}}^{c}$
$:=$ $X\backslash E_{\varepsilon _{1}}.$ Moreover, $|d\nu |$ $\ll $ $|d\mu
_{\varepsilon _{1}}|$ on $E_{\varepsilon _{1}}.$ Recall that from (\ref%
{eqn3.13}) we have%
\begin{equation}
d\mu _{\varepsilon _{j}}=((\varepsilon _{j}-\varepsilon _{1})A_{\varepsilon
_{1}}+N_{\varepsilon _{1}})|d\mu _{\varepsilon _{1}}|+(\varepsilon
_{j}-\varepsilon _{1})d\nu _{s}^{\varepsilon _{1}}  \label{eqn4.16}
\end{equation}

\noindent for $j$ $=$ $2,3,$ where $\varepsilon _{3}$ $\neq $ $\varepsilon
_{2}.$ By (\ref{eqn4.16}) we compute 
\begin{eqnarray}
&&\frac{|d\mu _{\varepsilon _{3}}|-|d\mu _{\varepsilon _{2}}|}{\varepsilon
_{3}-\varepsilon _{2}}  \label{4.16} \\
&=&\frac{|(\varepsilon _{3}-\varepsilon _{1})A_{\varepsilon
_{1}}+N_{\varepsilon _{1}}|-|(\varepsilon _{2}-\varepsilon
_{1})A_{\varepsilon _{1}}+N_{\varepsilon _{1}}|}{\varepsilon
_{3}-\varepsilon _{2}}|d\mu _{\varepsilon _{1}}|  \notag \\
&&+\frac{|\varepsilon _{3}-\varepsilon _{1}|-|\varepsilon _{2}-\varepsilon
_{1}|}{\varepsilon _{3}-\varepsilon _{2}}|d\nu _{s}^{\varepsilon _{1}}|. 
\notag
\end{eqnarray}%
\noindent Note that on $E_{\varepsilon _{1}},$ $N_{\varepsilon _{1}}$ $\neq $
$0$ a.e. [$|d\mu _{\varepsilon _{1}}|]$ and hence both $(\varepsilon
_{3}-\varepsilon _{1})A_{\varepsilon _{1}}+N_{\varepsilon _{1}}$ and $%
(\varepsilon _{2}-\varepsilon _{1})A_{\varepsilon _{1}}+N_{\varepsilon _{1}}$
cannot be zero simultaneously a.e. [$|d\mu _{\varepsilon _{1}}|]$ since $%
\varepsilon _{3}$ $\neq $ $\varepsilon _{2}.$ Therefore we can write%
\begin{eqnarray}
&&\frac{|(\varepsilon _{3}-\varepsilon _{1})A_{\varepsilon
_{1}}+N_{\varepsilon _{1}}|-|(\varepsilon _{2}-\varepsilon
_{1})A_{\varepsilon _{1}}+N_{\varepsilon _{1}}|}{\varepsilon
_{3}-\varepsilon _{2}}  \label{4.17} \\
&=&\frac{((\varepsilon _{3}-\varepsilon _{1})^{2}-(\varepsilon
_{2}-\varepsilon _{1})^{2})|A_{\varepsilon _{1}}|^{2}+2(\varepsilon
_{3}-\varepsilon _{2})A_{\varepsilon _{1}}\cdot N_{\varepsilon _{1}}}{%
(\varepsilon _{3}-\varepsilon _{2})(|(\varepsilon _{3}-\varepsilon
_{1})A_{\varepsilon _{1}}+N_{\varepsilon _{1}}|+|(\varepsilon
_{2}-\varepsilon _{1})A_{\varepsilon _{1}}+N_{\varepsilon _{1}}|)}  \notag \\
&=&\frac{((\varepsilon _{3}+\varepsilon _{2}-2\varepsilon
_{1})A_{\varepsilon _{1}}+2N_{\varepsilon _{1}})\cdot A_{\varepsilon _{1}}}{%
|(\varepsilon _{3}-\varepsilon _{1})A_{\varepsilon _{1}}+N_{\varepsilon
_{1}}|+|(\varepsilon _{2}-\varepsilon _{1})A_{\varepsilon
_{1}}+N_{\varepsilon _{1}}|}.  \notag
\end{eqnarray}

\noindent Suppose $(\varepsilon _{2}-\varepsilon _{1})A_{\varepsilon
_{1}}+N_{\varepsilon _{1}}$ $\neq $ $0$ on $E_{\varepsilon _{1}}$(recall
that $|d\mu _{\varepsilon _{1}}|$ is concentrated on $E_{\varepsilon _{1}}).$
Then it follows from (\ref{4.17}) that%
\begin{eqnarray}
&&\lim_{\varepsilon _{3}\rightarrow \varepsilon _{2}}\frac{|(\varepsilon
_{3}-\varepsilon _{1})A_{\varepsilon _{1}}+N_{\varepsilon
_{1}}|-|(\varepsilon _{2}-\varepsilon _{1})A_{\varepsilon
_{1}}+N_{\varepsilon _{1}}|}{\varepsilon _{3}-\varepsilon _{2}}  \label{4.18}
\\
&=&\frac{(\varepsilon _{2}-\varepsilon _{1})A_{\varepsilon
_{1}}+N_{\varepsilon _{1}}}{|(\varepsilon _{2}-\varepsilon
_{1})A_{\varepsilon _{1}}+N_{\varepsilon _{1}}|}\cdot A_{\varepsilon _{1}}. 
\notag
\end{eqnarray}

\noindent On $E_{\varepsilon _{1}},$ if $(\varepsilon _{2}-\varepsilon
_{1})A_{\varepsilon _{1}}+N_{\varepsilon _{1}}$ $=$ $0,$ then we have%
\begin{equation}
\frac{|d\mu _{\varepsilon _{3}}|-|d\mu _{\varepsilon _{2}}|}{\varepsilon
_{3}-\varepsilon _{2}}=\frac{|\varepsilon _{3}-\varepsilon _{2}|}{%
\varepsilon _{3}-\varepsilon _{2}}|A_{\varepsilon _{1}}||d\mu _{\varepsilon
_{1}}|  \label{4.19}
\end{equation}

\noindent by observing that $(\varepsilon _{3}-\varepsilon
_{1})A_{\varepsilon _{1}}+N_{\varepsilon _{1}}$ $=$ ($\varepsilon
_{3}-\varepsilon _{2})A_{\varepsilon _{1}}$ $+$ $(\varepsilon
_{2}-\varepsilon _{1})A_{\varepsilon _{1}}+N_{\varepsilon _{1}}$ $=$ ($%
\varepsilon _{3}-\varepsilon _{2})A_{\varepsilon _{1}}$ in (\ref{4.16}).
Also from (\ref{4.16}) we have%
\begin{equation}
\frac{|d\mu _{\varepsilon _{3}}|-|d\mu _{\varepsilon _{2}}|}{\varepsilon
_{3}-\varepsilon _{2}}=\frac{|\varepsilon _{3}-\varepsilon
_{1}|-|\varepsilon _{2}-\varepsilon _{1}|}{\varepsilon _{3}-\varepsilon _{2}}%
|d\nu _{s}^{\varepsilon _{1}}|.  \label{4.19.1}
\end{equation}

\noindent on $E_{\varepsilon _{1}}^{c}.$ Observe that 
\begin{eqnarray}
&&|\frac{|(\varepsilon _{3}-\varepsilon _{1})A_{\varepsilon
_{1}}+N_{\varepsilon _{1}}|-|(\varepsilon _{2}-\varepsilon
_{1})A_{\varepsilon _{1}}+N_{\varepsilon _{1}}|}{\varepsilon
_{3}-\varepsilon _{2}}|  \label{4.19.2} \\
&\leq &\frac{|((\varepsilon _{3}-\varepsilon _{1})A_{\varepsilon
_{1}}+N_{\varepsilon _{1}})-((\varepsilon _{2}-\varepsilon
_{1})A_{\varepsilon _{1}}+N_{\varepsilon _{1}})|}{|\varepsilon
_{3}-\varepsilon _{2}|}  \notag \\
&=&\frac{|(\varepsilon _{3}-\varepsilon _{2})A_{\varepsilon _{1}}|}{%
|\varepsilon _{3}-\varepsilon _{2}|}=|A_{\varepsilon _{1}}|  \notag
\end{eqnarray}%
\noindent by the triangle inequality. Since $d\mu $ and $d\nu $ are bounded
by assumption, we obtain that $|d\mu _{\varepsilon _{1}}|$ $=$ $|d\mu
+\varepsilon _{1}d\nu |$ is bounded. Note that $|A_{\varepsilon _{1}}|$ $\in 
$ $L^{1}(X,|d\mu _{\varepsilon _{1}}|)$ since \TEXTsymbol{\vert}$d\nu |$ $=$ 
$|A_{\varepsilon _{1}}||d\mu _{\varepsilon _{1}}|$ $+$ $|d\nu
_{s}^{\varepsilon _{1}}|$ from (\ref{eqn3.2}) and hence%
\begin{equation*}
\int_{X}|A_{\varepsilon _{1}}||d\mu _{\varepsilon _{1}}|\leq \int_{X}|d\nu
|<\infty
\end{equation*}
\noindent (note that $d\nu $ is bounded by assumption). Now from (\ref{4.16}%
), (\ref{4.18}), (\ref{4.19}), (\ref{4.19.1}), and (\ref{4.19.2}), we can
apply the Lebesgue dominated convergence theorem to conclude that%
\begin{eqnarray}
\mathcal{F}_{\pm }^{\prime }\mathcal{(}\varepsilon _{2})
&:&=\lim_{\varepsilon _{3}\rightarrow \varepsilon _{2}\pm }\int_{X}\frac{%
|d\mu _{\varepsilon _{3}}|-|d\mu _{\varepsilon _{2}}|}{\varepsilon
_{3}-\varepsilon _{2}}  \label{4.20} \\
&=&\int_{E_{\varepsilon _{1}}\cap \{(\varepsilon _{2}-\varepsilon
_{1})A_{\varepsilon _{1}}+N_{\varepsilon _{1}}\neq 0\}}\frac{(\varepsilon
_{2}-\varepsilon _{1})A_{\varepsilon _{1}}+N_{\varepsilon _{1}}}{%
|(\varepsilon _{2}-\varepsilon _{1})A_{\varepsilon _{1}}+N_{\varepsilon
_{1}}|}\cdot A_{\varepsilon _{1}}|d\mu _{\varepsilon _{1}}|  \notag \\
&&+\int_{E_{\varepsilon _{1}}\cap \{(\varepsilon _{2}-\varepsilon
_{1})A_{\varepsilon _{1}}+N_{\varepsilon _{1}}=0\}}\pm |A_{\varepsilon
_{1}}||d\mu _{\varepsilon _{1}}|  \notag \\
&&+\int_{E_{\varepsilon _{1}}^{c}}\frac{|\varepsilon _{2}-\varepsilon _{1}|}{%
\varepsilon _{2}-\varepsilon _{1}}|d\nu _{s}^{\varepsilon _{1}}|.  \notag
\end{eqnarray}

Comparing (\ref{4.20}) with (\ref{eqn1.5}) we obtain%
\begin{eqnarray}
&&\mathcal{F}_{\pm }^{\prime }\mathcal{(}\varepsilon _{2})-\mathcal{F}%
_{+}^{\prime }\mathcal{(}\varepsilon _{1})  \label{4.20.1} \\
&=&\int_{E_{\varepsilon _{1}}\cap \{(\varepsilon _{2}-\varepsilon
_{1})A_{\varepsilon _{1}}+N_{\varepsilon _{1}}\neq 0\}}\{(\frac{(\varepsilon
_{2}-\varepsilon _{1})A_{\varepsilon _{1}}+N_{\varepsilon _{1}}}{%
|(\varepsilon _{2}-\varepsilon _{1})A_{\varepsilon _{1}}+N_{\varepsilon
_{1}}|}-N_{\varepsilon _{1}})\cdot A_{\varepsilon _{1}}|d\mu _{\varepsilon
_{1}}|  \notag \\
&&+\int_{E_{\varepsilon _{1}}\cap \{(\varepsilon _{2}-\varepsilon
_{1})A_{\varepsilon _{1}}+N_{\varepsilon _{1}}=0\}}(\pm |A_{\varepsilon
_{1}}|-N_{\varepsilon _{1}}\cdot A_{\varepsilon _{1}})|d\mu _{\varepsilon
_{1}}|\}  \notag
\end{eqnarray}%
\noindent for $\varepsilon _{2}$ $>$ $\varepsilon _{1}$ (from which terms
involving $|d\nu _{s}^{\varepsilon _{1}}|$ cancel)$.$ On $E_{\varepsilon
_{1}}\cap \{(\varepsilon _{2}-\varepsilon _{1})A_{\varepsilon
_{1}}+N_{\varepsilon _{1}}\neq 0\}$ there holds%
\begin{equation}
|\frac{1}{\varepsilon _{2}-\varepsilon _{1}}(\frac{(\varepsilon
_{2}-\varepsilon _{1})A_{\varepsilon _{1}}+N_{\varepsilon _{1}}}{%
|(\varepsilon _{2}-\varepsilon _{1})A_{\varepsilon _{1}}+N_{\varepsilon
_{1}}|}-N_{\varepsilon _{1}})\cdot A_{\varepsilon _{1}}|\leq
2|A_{\varepsilon _{1}}|^{2}  \label{4.21}
\end{equation}%
\noindent by the same estimate as in deducing (\ref{eqn4.9}) (noting that $%
\varepsilon _{1}$ and $\varepsilon _{2}$ are not necessarily regular in the
estimate)$.$ On $E_{\varepsilon _{1}}\cap \{(\varepsilon _{2}-\varepsilon
_{1})A_{\varepsilon _{1}}+N_{\varepsilon _{1}}=0\}$ we have \TEXTsymbol{\vert%
}$\frac{1}{\varepsilon _{2}-\varepsilon _{1}}|$ $=$ $|A_{\varepsilon _{1}}|$
since \TEXTsymbol{\vert}$N_{\varepsilon _{1}}|$ $=$ $1.$ It follows that%
\begin{eqnarray}
&&\frac{\pm |A_{\varepsilon _{1}}|-N_{\varepsilon _{1}}\cdot A_{\varepsilon
_{1}}}{\varepsilon _{2}-\varepsilon _{1}}  \label{4.25} \\
&=&(\pm |A_{\varepsilon _{1}}|-N_{\varepsilon _{1}}\cdot A_{\varepsilon
_{1}})(sgn(\varepsilon _{2}-\varepsilon _{1})|A_{\varepsilon _{1}}|)  \notag
\\
&=&\pm sgn(\varepsilon _{2}-\varepsilon _{1})|A_{\varepsilon
_{1}}|^{2}+|A_{\varepsilon _{1}}|^{2}  \notag \\
&=&2|A_{\varepsilon _{1}}|^{2}\text{ \ if }\varepsilon _{2}>(<\text{, resp.})%
\text{ }\varepsilon _{1}\text{ for the case of + (}-\text{, resp.) sign} 
\notag \\
( &=&0\text{ if }\varepsilon _{2}>(<\text{, resp.})\text{ }\varepsilon _{1}%
\text{ for the case of }-\text{(}+\text{, resp.) sign).}  \notag
\end{eqnarray}%
\noindent Here we have used the fact that $N_{\varepsilon _{1}}\cdot
A_{\varepsilon _{1}}$ $=$ $-(\varepsilon _{2}$ $-$ $\varepsilon _{1})^{-1}$ $%
=$ $sgn(\varepsilon _{1}$ $-$ $\varepsilon _{2})|A_{\varepsilon _{1}}|$ on $%
E_{\varepsilon _{1}}\cap \{(\varepsilon _{2}-\varepsilon _{1})A_{\varepsilon
_{1}}+N_{\varepsilon _{1}}=0\}$ (in which \TEXTsymbol{\vert}$\frac{1}{%
\varepsilon _{2}-\varepsilon _{1}}|$ $=$ $|A_{\varepsilon _{1}}|).$ From (%
\ref{4.20.1}) and (\ref{4.25}), we can write%
\begin{eqnarray}
&&\frac{\mathcal{F}_{\pm }^{\prime }\mathcal{(}\varepsilon _{2})-\mathcal{F}%
_{+}^{\prime }\mathcal{(}\varepsilon _{1})}{\varepsilon _{2}-\varepsilon _{1}%
}  \label{4.26} \\
&=&\int_{E_{\varepsilon _{1}}\cap \{(\varepsilon _{2}-\varepsilon
_{1})A_{\varepsilon _{1}}+N_{\varepsilon _{1}}\neq 0\}}g_{(\varepsilon
_{1},\varepsilon _{2})}|d\mu _{\varepsilon _{1}}|+\int_{E_{\varepsilon
_{1}}\cap \{(\varepsilon _{2}-\varepsilon _{1})A_{\varepsilon
_{1}}+N_{\varepsilon _{1}}=0\}}h_{(\varepsilon _{1})}|d\mu _{\varepsilon
_{1}}|  \notag
\end{eqnarray}%
\noindent where 
\begin{equation}
g_{(\varepsilon _{1},\varepsilon _{2})}=\frac{1}{\varepsilon
_{2}-\varepsilon _{1}}(\frac{(\varepsilon _{2}-\varepsilon
_{1})A_{\varepsilon _{1}}+N_{\varepsilon _{1}}}{|(\varepsilon
_{2}-\varepsilon _{1})A_{\varepsilon _{1}}+N_{\varepsilon _{1}}|}%
-N_{\varepsilon _{1}})\cdot A_{\varepsilon _{1}},  \label{4.27}
\end{equation}%
\noindent and%
\begin{eqnarray}
&&h_{(\varepsilon _{1})}  \label{4.28} \\
&=&2|A_{\varepsilon _{1}}|^{2}\text{ if }\varepsilon _{2}>(<\text{, resp.})%
\text{ }\varepsilon _{1}\text{ in the case of + (}-\text{, resp.) sign} 
\notag \\
( &=&0\text{ if }\varepsilon _{2}>(<\text{, resp.})\text{ }\varepsilon _{1}%
\text{ in the case of }-\text{(}+\text{, resp.) sign).}  \notag
\end{eqnarray}%
\noindent From (\ref{4.21}) and (\ref{4.28}), we have%
\begin{eqnarray}
|g_{(\varepsilon _{1},\varepsilon _{2})}| &\leq &2|A_{\varepsilon _{1}}|^{2}%
\text{ on }E_{\varepsilon _{1}}\cap \{(\varepsilon _{2}-\varepsilon
_{1})A_{\varepsilon _{1}}+N_{\varepsilon _{1}}\neq 0\},  \label{4.29} \\
|h_{(\varepsilon _{1})}| &\leq &2|A_{\varepsilon _{1}}|^{2}\text{ on }%
E_{\varepsilon _{1}}\cap \{(\varepsilon _{2}-\varepsilon _{1})A_{\varepsilon
_{1}}+N_{\varepsilon _{1}}=0\}.  \notag
\end{eqnarray}%
\noindent Now given a point $p$ $\in $ $E_{\varepsilon _{1}}$ (modulo a $%
|d\mu _{\varepsilon _{1}}|$-measure zero set)$,$ there is at most one $%
\varepsilon _{2}$ $\neq $ $\varepsilon _{1}$ such that $(\varepsilon
_{2}-\varepsilon _{1})A_{\varepsilon _{1}}(p)$ $+$ $N_{\varepsilon _{1}}(p)$ 
$=$ $0.$ The reason is that if there are two distinct such $\varepsilon
_{2}, $ then $A_{\varepsilon _{1}}(p)$ $=$ $0$ and hence $N_{\varepsilon
_{1}}(p)$ $=$ $0.$ So all such points form a $|d\mu _{\varepsilon _{1}}|$%
-measure zero set since $N_{\varepsilon _{1}}$ $\neq $ $0$ a.e. $[|d\mu
_{\varepsilon _{1}}|].$ We denote such $\varepsilon _{2}$ by $\varepsilon
_{2}(p).$ Then for any $\varepsilon ,$ $\varepsilon _{1}$ $<$ $\varepsilon $ 
$<$ $\varepsilon _{2}(p),$ there holds $(\varepsilon -\varepsilon
_{1})A_{\varepsilon _{1}}(p)+N_{\varepsilon _{1}}(p)\neq 0.$ So we have%
\begin{eqnarray}
\lim_{\varepsilon \rightarrow \varepsilon _{1}+}\tilde{g}_{(\varepsilon
_{1},\varepsilon )}(p) &=&\lim_{\varepsilon \rightarrow \varepsilon
_{1}+}g_{(\varepsilon _{1},\varepsilon )}(p)  \label{4.30} \\
&=&f^{\prime }(0+)\cdot A_{\varepsilon _{1}}(p)  \notag \\
&=&|A_{\varepsilon _{1}}(p)|^{2}-|(A_{\varepsilon _{1}}\cdot N_{\varepsilon
_{1}})(p)|^{2}  \notag
\end{eqnarray}%
\noindent by (\ref{eqn4.11}), where $\tilde{g}_{(\varepsilon
_{1},\varepsilon )}$ is defined on $E_{\varepsilon _{1}}$ as follows:%
\begin{eqnarray*}
\tilde{g}_{(\varepsilon _{1},\varepsilon )} &=&g_{(\varepsilon
_{1},\varepsilon )}\text{ on }E_{\varepsilon _{1}}\cap \{(\varepsilon
-\varepsilon _{1})A_{\varepsilon _{1}}+N_{\varepsilon _{1}}\neq 0\}\text{ and%
} \\
&=&h_{(\varepsilon _{1})}\text{ on }E_{\varepsilon _{1}}\cap \{(\varepsilon
-\varepsilon _{1})A_{\varepsilon _{1}}+N_{\varepsilon _{1}}=0\}.
\end{eqnarray*}%
\noindent In view of (\ref{4.29}), (\ref{4.30}), and the assumption $%
|A_{\varepsilon _{1}}|^{2}$ $\in $ $L^{1}(X,|d\mu _{\varepsilon _{1}}|)$, we
can now apply the Lebesgue dominated convergence theorem to compute%
\begin{eqnarray}
&&\lim_{\varepsilon \rightarrow \varepsilon _{1}+}\frac{\mathcal{F}_{\pm
}^{\prime }\mathcal{(}\varepsilon )-\mathcal{F}_{+}^{\prime }\mathcal{(}%
\varepsilon _{1})}{\varepsilon -\varepsilon _{1}}  \label{4.31} \\
&=&\lim_{\varepsilon \rightarrow \varepsilon _{1}+}\int_{E_{\varepsilon
_{1}}}\tilde{g}_{(\varepsilon _{1},\varepsilon )}|d\mu _{\varepsilon _{1}}| 
\notag \\
&=&\int_{E_{\varepsilon _{1}}}(\lim_{\varepsilon \rightarrow \varepsilon
_{1}+}\tilde{g}_{(\varepsilon _{1},\varepsilon )})|d\mu _{\varepsilon _{1}}|
\notag \\
&=&\int_{X}(|A_{\varepsilon _{1}}|^{2}-|(A_{\varepsilon _{1}}\cdot
N_{\varepsilon _{1}})|^{2})|d\mu _{\varepsilon _{1}}|\geq 0.  \notag
\end{eqnarray}

\noindent In the last equality of (\ref{4.31}), we have used the fact that $%
|d\mu _{\varepsilon _{1}}|(X\backslash E_{\varepsilon _{1}})$ $=$ $0$ since $%
|d\mu _{\varepsilon _{1}}|$ is concentrated on $E_{\varepsilon _{1}}.$
Similarly we also have

\begin{equation*}
\lim_{\varepsilon \rightarrow \varepsilon _{1}-}\frac{\mathcal{F}_{\pm
}^{\prime }\mathcal{(}\varepsilon )-\mathcal{F}_{-}^{\prime }\mathcal{(}%
\varepsilon _{1})}{\varepsilon -\varepsilon _{1}}=\int_{X}\{|A_{\varepsilon
_{1}}|^{2}-|(A_{\varepsilon _{1}}\cdot N_{\varepsilon _{1}})|^{2}\}|d\mu
_{\varepsilon _{1}}|\geq 0.
\end{equation*}

\noindent Note that $\frac{|\varepsilon _{2}-\varepsilon _{1}|}{\varepsilon
_{2}-\varepsilon _{1}}|d\nu _{s}^{\varepsilon _{1}}|$ in (\ref{4.20})
cancels in both cases. We have proved Theorem C (2) (\ref{eqn1.7}).

\bigskip

We remark that for $\varepsilon _{2}$ $\neq $ $\varepsilon _{1},$ $|d\nu
_{s}^{\varepsilon _{2}}|$ is concentrated on $E_{\varepsilon _{1}}\cap
\{(\varepsilon _{2}-\varepsilon _{1})A_{\varepsilon _{1}}+N_{\varepsilon
_{1}}=0\}.$ Since $|d\nu _{s}^{\varepsilon _{2}}|$ $\perp $ $|d\nu
_{s}^{\varepsilon _{1}}|$ by Lemma 3.1 and $|d\mu _{\varepsilon _{1}}|$ $%
\perp $ $|d\nu _{s}^{\varepsilon _{1}}|,$ we obtain that $|d\nu
_{s}^{\varepsilon _{2}}|$ is concentrated on $E_{\varepsilon _{1}}$ by Lemma
3.2 (1) or (\ref{3.7}). Moreover, observing that $|d\nu _{s}^{\varepsilon
_{2}}|$ $\perp $ $|d\mu _{\varepsilon _{2}}|,$ we conclude that $|d\nu
_{s}^{\varepsilon _{2}}|$ is concentrated on $E_{\varepsilon _{1}}\cap
\{(\varepsilon _{2}-\varepsilon _{1})A_{\varepsilon _{1}}+N_{\varepsilon
_{1}}=0\}$ in view of (\ref{eqn3.13}).

\bigskip

\textbf{Example 4.1. }Continue the discussion in Example 3.1. Since every $%
\varepsilon $ is regular in this case, by Theorem C (1) and (\ref{eqn3.21})
we have%
\begin{eqnarray*}
\mathcal{F}^{\prime \prime }(0) &=&\int_{X}\{|A_{0}|^{2}-|(A_{0}\cdot
N_{0})|^{2}\}|d\mu | \\
&=&\int_{X}\frac{|\nabla v|^{2}+(|\nabla v|^{2}|\nabla u|^{2}-|(\nabla
v\cdot \nabla u)|^{2})}{(1+|\nabla u|^{2})^{3/2}}d^{m}x\geq 0
\end{eqnarray*}

\noindent by Cauchy's inequality. This implies that a Riemannian minimal
graph in $R^{m+1}$ over $X$ $\subset $ $R^{m}$ has the local area-minimizing
property.

\bigskip

\textbf{Example 4.2.} Continue the discussion in Example 3.2. Suppose that $%
u $ and $\varphi $ are in $C^{1,1}$ or $C^{2}.$ Then the singular set of the
graph defined by $u$ $+$ $\varepsilon \varphi $ has vanishing Lebesgue
measure in $R^{2n}$ according to \cite{Ba}. It follows that $d\nu
_{s}^{\varepsilon }$ $=$ $0,$ and hence each $\varepsilon $ is regular in
this situation. By Theorem C (1) and (\ref{eqn3.23}) we have

\begin{equation*}
\mathcal{F}^{\prime \prime }(0)=\int_{\Omega \backslash S(u)}\frac{|\nabla u-%
\vec{X}^{\ast }|^{2}|\nabla \varphi |^{2}-|(\nabla u-\vec{X}^{\ast })\cdot
\nabla \varphi |^{2}}{|\nabla u-\vec{X}^{\ast }|^{3}}d^{2n}x\geq 0
\end{equation*}

\noindent by Cauchy's inequality again. So a $C^{1,1}$ or $C^{2}$ $p$-area
stationary graph in the Heisenberg group over $\Omega $ $\subset $ $R^{2n}$
has the local $p$- area- minimizing property. This fact was shown by a
calibration argument in \cite{chmy} for the nonsingular case with $n=1$.
Later Ritor\'{e} and Rosales (\cite{Ri2}) extended the result to the
situation having singularities.

\bigskip

\textbf{Example 4.3. }Continue the discussion in Example 3.3. In \cite{chmy}
we computed the second variation of the $p$-area in the direction of a
vector field with support away from the singular set. Here we consider the
situation of variations with support containing a singular curve as in
Example 3.3. The following computation is based on a private talk given by
Hung-Lin Chiu. Recall that $v$ $=$ $v_{2}e_{2}$ $+$ $fT$ (take $v_{1}$ $=$ $%
0 $ for simplicity). Then we follow the argument in \cite{chmy} to get

\begin{eqnarray}
\delta _{v}^{2}\int_{\Sigma }\Theta \wedge e^{1} &=&\int_{\Sigma \backslash
S_{\Sigma }}L_{v}^{2}(\Theta \wedge e^{1})  \label{eqn4.23} \\
&=&\int_{\Sigma \backslash S_{\Sigma }}L_{v}\{(f\alpha
-v_{2})H+d(fe^{1})\}\Theta \wedge e^{1}  \notag \\
&=&\int_{\Sigma \backslash S_{\Sigma }}-(f\alpha -v_{2})^{2}(e_{2}H)\Theta
\wedge e^{1}+d\circ L_{v}(fe^{1})  \notag
\end{eqnarray}

\noindent for a $p$-area stationary surface $\Sigma $ (hence $H$ $=$ $0$ on $%
\Sigma ).$ We can express the first term of the last integrand in (\ref%
{eqn4.23}) in terms of pseudohermitian geometric quantities (see Section 6
in \cite{chmy}). The second term of the same integrand reflects the
contribution of the singular curve $S_{\Sigma }$ as shown below. By a direct
computation we obtain%
\begin{eqnarray}
L_{v}(fe^{1}) &=&(vf+f^{2}\func{Re}A^{1}\text{ }_{\bar{1}})e^{1}+(\omega (T)+%
\func{Im}A^{1}\text{ }_{\bar{1}})f^{2}e^{2}  \label{eqn4.24} \\
&&-(\omega (T)+\func{Im}A^{1}\text{ }_{\bar{1}})fv_{2}\Theta  \notag
\end{eqnarray}

\noindent where $A^{1}$ $_{\bar{1}}$ and $\omega $ denote the
pseudohermitian torsion and connection form, respectively. Since $\Theta
(\tau )$ $=$ $0$, $e_{+}^{1}(\tau )$ $+$ $e_{-}^{1}(-\tau )$ $=$ $0,$ and $%
e_{+}^{2}(\tau )=-e_{-}^{2}(\tau )$ by (\ref{eqn3.30}), from (\ref{eqn4.24})
we have%
\begin{eqnarray}
\int_{\Sigma \backslash S_{\Sigma }}d\circ L_{v}(fe^{1}) &=&\int_{S_{\Sigma
}^{+}}(\omega (T)+\func{Im}A^{1}\text{ }_{\bar{1}})f^{2}e_{+}^{2}
\label{eqn4.25} \\
&&+\int_{S_{\Sigma }^{-}}(\omega (T)+\func{Im}A^{1}\text{ }_{\bar{1}%
})f^{2}e_{-}^{2}  \notag \\
&=&2\int_{S_{\Sigma }^{+}}(\omega (T)+\func{Im}A^{1}\text{ }_{\bar{1}%
})f^{2}e_{+}^{2}(\tau )ds  \notag
\end{eqnarray}

\noindent where $s$ is the unit-speed parameter for $S_{\Sigma }^{+}.$

\bigskip

\section{\textbf{Appendix: Generalized Heisenberg geometry}}

We will discuss the notions of gradient and hypersurface area in a general
formulation unifying Riemannian and pseudohermitian (horizontal or
Heisenberg) structures (see, e.g., \cite{St}).

Let $M$ be an $m-$dimensional differentiable manifold with a nonnegative
inner product $<\cdot ,\cdot >$ on its cotangent bundle $T^{\ast }M.$
Namely, $<\cdot ,\cdot >$ is a symmetric bilinear form such that $<\omega
,\omega >$ $\geq $ $0$ for any $\omega $ $\in $ $T^{\ast }M.$ Some authors
call such a manifold $M$ subriemannian. Clearly if $<\cdot ,\cdot >$ is
positive definite, $(M,<\cdot ,\cdot >)$ is a Riemannian manifold. For $M$
being the Heisenberg group $H_{n}$ of dimension $m$ $=$ $2n+1,$ let 
\begin{equation*}
\hat{e}_{j}=\frac{\partial }{\partial x_{j}}+y_{j}\frac{\partial }{\partial z%
},\hat{e}_{j^{\prime }}=\frac{\partial }{\partial y_{j}}-x_{j}\frac{\partial 
}{\partial z},
\end{equation*}

\noindent $1\leq j\leq n$ be the left-invariant vector fields on $H_{n},$ in
which $x_{1},$ $y_{1},$ $x_{2},$ $y_{2},$ $...,$ $x_{n},$ $y_{n},$ $z$
denote the coordinates of $H_{n}.$ The (contact) 1-form $\Theta $ $\equiv $ $%
dz+\sum_{j=1}^{n}(x_{j}dy_{j}-y_{j}dx_{j})$ annihilates $\hat{e}_{j}^{\prime
}s$ and $\hat{e}_{j^{\prime }}^{\prime }s.$ We observe that $dx_{1},$ $%
dy_{1},$ $dx_{2},$ $dy_{2},$ $...,$ $dx_{n},$ $dy_{n},$ $\Theta $ are dual
to $\hat{e}_{1},$ $\hat{e}_{1^{\prime }},$ $\hat{e}_{2},$ $\hat{e}%
_{2^{\prime }},$ ...,$\hat{e}_{n},$ $\hat{e}_{n^{\prime }},$ $\frac{\partial 
}{\partial z}.$ Define a nonnegative inner product by%
\begin{eqnarray}
&<&dx_{j},dx_{k}>=\delta _{jk},\text{ }<dy_{j},dy_{k}>=\delta
_{jk},<dx_{j},dy_{k}>=0,  \label{H1} \\
&<&\Theta ,dx_{j}>=<\Theta ,dy_{k}>=<\Theta ,\Theta >=0.  \notag
\end{eqnarray}

We can extend the definition of the above nonnegative inner product to the
situation of a general pseudohermitian manifold. Take $e_{j},$ $e_{j^{\prime
}}=Je_{j}$, $j=1,$ $2,$ $...,n$ to be an orthonormal basis in the kernel of
the contact form $\Theta $ with respect to the Levi metric $\frac{1}{2}%
d\Theta (\cdot ,J\cdot ).$ Let $T$ be the Reeb vector field of $\Theta $
(such that $\Theta (T)$ $=$ $1$ and $d\Theta (T,\cdot )$ $=$ $0).$ Denote
the dual coframe of $e_{j},$ $e_{j^{\prime }},$ $T$ by $\theta ^{j},$ $%
\theta ^{j^{\prime }}$ (and $\Theta ).$ Now we can replace $dx_{j}$, $dy_{j}$
by $\theta ^{j},$ $\theta ^{j^{\prime }}$ in (\ref{H1}) to define a
nonnegative inner product on a general pseudohermitian manifold:

\begin{eqnarray}
&<&\theta ^{j},\theta ^{k}>=\delta _{jk},\text{ }<\theta ^{j^{\prime
}},\theta ^{k^{\prime }}>=\delta _{jk},<\theta ^{j},\theta ^{k^{\prime }}>=0,
\label{H2} \\
&<&\Theta ,\theta ^{j}>=<\Theta ,\theta ^{k^{\prime }}>=<\Theta ,\Theta >=0.
\notag
\end{eqnarray}

We use the same notation $<\cdot ,\cdot >$ to denote the pairing between $TM$
and $T^{\ast }M.$ Define the bundle morphism $G:T^{\ast }M\rightarrow TM$ by 
\begin{equation}
<G(\omega ),\eta >=<\omega ,\eta >  \label{H3}
\end{equation}

\noindent for $\omega ,\eta $ $\in $ $T^{\ast }M.$ In the Riemannian case, $%
G $ is in fact an isometry. In the pseudohermitian case, $G(T^{\ast }M)$ is
the contact subbundle $\xi $ of $TM,$ the kernel of $\Theta .$ By letting $%
\eta $ $=$ $\Theta $ in (\ref{H3})$,$ we get $G(T^{\ast }M)$ $\subset $ $\xi
.$ On the other hand, it is easy to see that $G(\theta ^{j})$ $=$ $e_{j},$ $%
G(\theta ^{j^{\prime }})$ $=$ $e_{j^{\prime }}$ (and $G(\Theta )$ $=$ $0)$.
Since $e_{j},$ $e_{j^{\prime }}$, $j=1,$ $2,$ $...,n$ span $\xi ,$ we have $%
\xi $ $\subset $ $G(T^{\ast }M).$ For a smooth function $\varphi $ on $M,$
we define the gradient $\nabla \varphi $ $:=$ $G(d\varphi ).$ In the
pseudohermitian case, this $\nabla \varphi $ is nothing but the subgradient $%
\nabla _{b}\varphi $ $:=$ $\sum_{j=1}^{n}\{e_{j}(\varphi )e_{j}$ $+$ $%
e_{j^{\prime }}(\varphi )e_{j^{\prime }}\}.$

Let $M$ be a general subriemannian manifold of dimension $n+1$, i.e., an ($%
n+1)$-dimensional differentiable manifold with a nonnegative inner product $%
<\cdot ,\cdot >$ on its cotangent bundle $T^{\ast }M.$ Let $\varphi $ be a
defining function of a hypersurface $\Sigma $ $\subset $ $M.$ That is, $%
\Sigma $ $=$ $\{\varphi =0\}.$ Given a volume form $dv_{M}$ (independent of $%
<\cdot ,\cdot >$), we can define an area (or volume) element $dv_{\Sigma }$
of $\Sigma $ up to sign by%
\begin{equation}
dv_{\Sigma }=\frac{d\varphi }{|d\varphi |}\text{ }\rfloor \text{ }dv_{M}
\label{H4}
\end{equation}

\noindent restricted to $\Sigma .$ Here for $\omega ,$ $\eta $ $\in $ $%
T^{\ast }M,$ $|\omega |$ $:=$ $<\omega ,\omega >^{1/2}$ and $\omega $ $%
\rfloor $ $dv_{M}$ is defined so that 
\begin{equation}
\eta \wedge (\omega \rfloor dv_{M})=<\eta ,\omega >dv_{M}.  \label{H4.1}
\end{equation}%
\noindent If we write $dv_{M}$ $=$ $\omega ^{1}$ $\wedge $ $\omega ^{2}$ $%
\wedge $ $...$ $\wedge $ $\omega ^{n+1}$ for independent $1$-forms $\omega
^{j}$'s and $\omega $ $=$ $\lambda _{j}\omega ^{j}$ (summation convention)$,$
then it is straightforward to verify that%
\begin{eqnarray}
&&\omega \rfloor dv_{M}  \label{H4.1.5} \\
&=&\lambda _{j}<\omega ^{j},\omega ^{k}>(-1)^{k-1}\omega ^{1}\wedge ..\wedge 
\hat{\omega}^{k}\wedge ..\wedge \omega ^{n+1}  \notag
\end{eqnarray}%
\noindent satisfies $\omega ^{l}\wedge (\omega \rfloor dv_{M})$ $=$ $<\omega
^{l},\omega >dv_{M}$ for all $l$ and hence (\ref{H4.1}) holds for all $\eta
. $ On the other hand, there is a unique $n$-form $\Phi $ satisfying 
\begin{equation}
\eta \wedge \Phi =<\eta ,\omega >dv_{M}  \label{H4.2}
\end{equation}%
\noindent for all $1$-forms $\eta .$ Suppose there are two $n$-forms $\Phi
_{1},$ $\Phi _{2}$ satisfying (\ref{H4.2})$.$ Then it follows that $\eta
\wedge (\Phi _{1}-\Phi _{2})$ $=$ $0$ for all $1$-forms $\eta $ and hence $%
\Phi _{1}-\Phi _{2}$ $=$ $0$. We have justified the formula (\ref{H4.1.5}).
There is an intrinsic expression for $\omega \rfloor dv_{M}$ as follows:%
\begin{eqnarray*}
&&(\omega \rfloor dv_{M})(X_{1},...,X_{n}) \\
&=&dv_{M}(G(\omega ),X_{1},...,X_{n}).
\end{eqnarray*}%
\noindent Note that $dv_{\Sigma }$ defined by (\ref{H4}) is independent of
the choice of $\varphi $ by a positive scalar multiple function, but changes
sign if $\varphi $ is replaced by $-\varphi $. Now we write $dv_{M}$ $=$ $%
\omega ^{1}$ $\wedge $ $\omega ^{2}$ $\wedge $ $...$ $\wedge $ $\omega
^{n+1} $ for independent $1$-forms $\omega ^{j}$'s and compute

\begin{eqnarray}
d\varphi \rfloor dv_{M} &=&d\varphi \rfloor \omega ^{1}\wedge \omega
^{2}\wedge ...\wedge \omega ^{n+1}  \label{H5} \\
&=&v_{i}(\varphi )<\omega ^{i},\omega ^{j}>(-1)^{j-1}\omega ^{1}\wedge ..%
\hat{\omega}^{j}..\wedge \omega ^{n+1}  \notag
\end{eqnarray}%
\noindent by (\ref{H4.1.5}), where $v_{i}$'s are tangent vectors dual to $%
\omega ^{j}$'s and $\hat{\omega}^{j}$ means $\omega ^{j}$ deleted. On $%
\Sigma ,$ $d\varphi $ $=$ $v_{i}(\varphi )\omega ^{i}$ $=$ $0.$ Assuming $%
v_{n+1}(\varphi )$ $\neq $ $0,$ say, we have

\begin{equation}
\omega ^{n+1}=-\frac{1}{v_{n+1}(\varphi )}\sum_{j=1}^{n}v_{j}(\varphi
)\omega ^{j}.  \label{H6}
\end{equation}%
\noindent Substituting (\ref{H6}) into (\ref{H5}) and noting that $|d\varphi
|^{2}$ $=$ $v_{i}(\varphi )<\omega ^{i},\omega ^{j}>v_{j}(\varphi ),$ we
obtain

\begin{equation}
dv_{\Sigma }=\frac{(-1)^{n}}{v_{n+1}(\varphi )}|d\varphi |\omega ^{1}\wedge
\omega ^{2}\wedge ...\wedge \omega ^{n}.  \label{H7}
\end{equation}

\bigskip

\textbf{Example A.1.} Suppose $\Sigma $ is a hypersurface of $M$ $=$ $%
R^{n+1}.$ Take $<\cdot ,\cdot >$ and $dv_{M}$ to be the Euclidean metric and
the associated volume form, respectively. Write the defining function $%
\varphi $ $=$ $z-u(x_{1},$ $x_{2},$ $...,$ $x_{n})$ and $dv_{M}$ $=$ $dx_{1}$
$\wedge $ $dx_{2}$ $\wedge $ $...$ $\wedge $ $dx_{n}$ $\wedge $ $dz$ where $%
x_{1},$ $x_{2},$ $...,$ $x_{n},$ $z$ are coordinates of $R^{n+1}.$ It
follows that $|d\varphi |$ $=$ $(1$ $+$ $u_{x^{1}}^{2}$ $+$ $...$ $+$ $%
u_{x^{n}}^{2})^{1/2}$ and $d\varphi $ $\rfloor $ $dx_{1}$ $\wedge $ $dx_{2}$ 
$\wedge $ $...$ $\wedge $ $dx_{n}$ $\wedge $ $dz$ $=$ $(-1)^{n}(1$ $+$ $%
u_{x^{1}}^{2}$ $+$ $...$ $+$ $u_{x^{n}}^{2})$ $dx_{1}$ $\wedge $ $dx_{2}$ $%
\wedge $ $...$ $\wedge $ $dx_{n}$ when restricted to $\Sigma .$ So taking $%
\omega ^{j}$ $=$ $dx_{j}$, $j$ $=$ $1,$ $...$ $,$ $n,$ $\omega ^{n+1}$ $=$ $%
dz$ and noting that $v_{n+1}$ $=$ $\frac{\partial }{\partial z},$ $%
v_{n+1}(\varphi )$ $=$ $1$ in (\ref{H7})$,$ we have

\begin{equation*}
dv_{\Sigma }=(-1)^{n}(1+u_{x^{1}}^{2}+...+u_{x^{n}}^{2})^{1/2}dx_{1}\wedge
dx_{2}\wedge ...\wedge dx_{n}.
\end{equation*}

\noindent This is the standard area element (up to a sign) for a graph in
Euclidean space.

\bigskip

\textbf{Example A.2. }For $M$ being the Heisenberg group of dimension $2n+1$%
, we take the volume form $dv_{M}$ $=$ $dx_{1}$ $\wedge $ $dy_{1}$ $\wedge $ 
$...\wedge $ $dx_{n}$ $\wedge $ $dy_{n}$ $\wedge $ $\Theta $ (the volume
form with respect to the left invariant metric)$.$ Let $\omega ^{2j-1}$ $=$ $%
dx_{j},$ $\omega ^{2j}$ $=$ $dy_{j}$, $1$ $\leq $ $j$ $\leq $ $n,$ $\omega
^{n+1}$ $=$ $\Theta $ while $v_{2j-1}$ $=$ $\hat{e}_{j},$ $v_{2j}$ $=$ $\hat{%
e}_{j^{\prime }},$and $v_{2n+1}$ $=$ $\frac{\partial }{\partial z}.$ For $%
\varphi $ $=$ $z-u(x_{1},$ $y_{1},$ $...,$ $x_{n},$ $y_{n})$, we compute $%
v_{2n+1}(\varphi )$ $=$ $1$ and $|d\varphi |^{2}$ $=$ $\sum_{j=1}^{n}\{\hat{e%
}_{j}(\varphi )^{2}$ $+$ $\hat{e}_{j^{\prime }}(\varphi )^{2}\}$ $=$ $%
\sum_{j=1}^{n}\{(u_{x_{j}}-y_{j})^{2}$ $+$ $(u_{y_{j}}+x_{j})^{2}\}$ by (\ref%
{H1})$.$ Substituting these formulas into (\ref{H7}) gives%
\begin{equation*}
dv_{\Sigma
}=[\sum_{j=1}^{n}\{(u_{x_{j}}-y_{j})^{2}+(u_{y_{j}}+x_{j})^{2}%
\}]^{1/2}dx_{1}\wedge dy_{1}\wedge ...\wedge dx_{n}\wedge dy_{n}
\end{equation*}

\noindent (note that $(-1)^{2n}$ $=$ $1).$ This is the standard ($p$- or $H$%
-) area element for a graph in the Heisenberg group.

We can also recover the area element of an intrinsic graph (e.g., \cite{ASCV}%
, \cite{BASCV}) in the Heisenberg group from (\ref{H7}). Let us explain this
for the 3-dimensional case ($n=2)$.

\bigskip

\textbf{Example A.3. }Take $\omega ^{1}$ $=$ $dy,$ $\omega ^{2}$ $=$ $\Theta 
$ $=$ $dz$ $+$ $xdy$ $-$ $ydx,$ and $\omega ^{3}$ $=$ $dx$ ($x$ $=$ $x_{1},$ 
$y$ $=$ $y_{1})$ (so $v_{1}$ $=$ $\hat{e}_{1^{\prime }},$ $v_{2}$ $=$ $\frac{%
\partial }{\partial z},$ and $v_{3}$ $=$ $\hat{e}_{1}).$ We compute $%
|d\varphi |^{2}$ $=$ $\hat{e}_{1}(\varphi )^{2}$ $+$ $\hat{e}_{1^{\prime
}}(\varphi )^{2}$ by (\ref{H1}) and reduce (\ref{H7}) to%
\begin{equation}
dv_{\Sigma }=\sqrt{1+(\frac{\hat{e}_{1^{\prime }}(\varphi )}{\hat{e}%
_{1}(\varphi )})^{2}}dy\wedge \Theta .  \label{H8}
\end{equation}

\noindent An intrinsic graph is parametrized by $\eta ,$ $\tau $ as follows:
(we have adjusted the normalization constant)%
\begin{equation}
x=\phi (\eta ,\tau ),y=\eta ,z=\tau +\eta \phi (\eta ,\tau )  \label{H9}
\end{equation}

\noindent It follows that $\Theta $ $=$ $d\tau $ $+$ $2\phi d\eta ,$ $dy$ $=$
$d\eta ,$ and hence $dy$ $\wedge $ $\Theta $ $=$ $d\eta $ $\wedge $ $d\tau $
by (\ref{H9})$.$ In coordinates $(\rho ,$ $\eta ,$ $\tau )$ related to $(x,$ 
$y,$ $z)$ by $x$ $=$ $\rho ,$ $y$ $=$ $\eta ,$ $z$ $=$ $\tau $ $+\ \eta \rho
,$ we can write the defining function $\varphi $ $=$ $\rho $ $-$ $\phi (\eta
,\tau ).$ By the chain rule we obtain $\frac{\partial }{\partial x}$ $=$ $%
\frac{\partial }{\partial \rho }$ $-$ $\eta \frac{\partial }{\partial \tau }%
, $ $\frac{\partial }{\partial y}$ $=$ $\frac{\partial }{\partial \eta }$ $-$
$\rho \frac{\partial }{\partial \tau },$ and $\frac{\partial }{\partial z}$ $%
= $ $\frac{\partial }{\partial \tau }.$ It follows that $\hat{e}_{1}$ $=$ $%
\frac{\partial }{\partial \rho },$ $\hat{e}_{1^{\prime }}$ $=$ $\frac{%
\partial }{\partial \eta }$ $-$ $2\rho \frac{\partial }{\partial \tau },$
and hence $\hat{e}_{1}(\varphi )$ $=$ $1$, $\hat{e}_{1^{\prime }}(\varphi )$ 
$=$ $-\phi _{\eta }$ $+$ $2\rho \phi _{\tau }.$ Substituting these formulas
into (\ref{H8})(and noting that $\rho $ $=$ $\phi (\eta ,\tau )$ when
restricted to $\Sigma )$, we obtain%
\begin{equation}
dv_{\Sigma }=\sqrt{1+(\phi _{\eta }-2\phi \phi _{\tau })^{2}}d\eta \wedge
d\tau  \label{H10}
\end{equation}

\noindent (e.g., \cite{ASCV}, \cite{BASCV}).

\bigskip

Next we consider $\omega ^{j}$'s to be a moving coframe such that $%
v_{j}(\varphi )$ $=$ $0$ for $1$ $\leq $ $j$ $\leq $ $n,$ and $%
v_{n+1}(\varphi )$ $\neq $ $0$ in (\ref{H7})$.$ It follows that $|d\varphi |$
$=$ $|v_{n+1}(\varphi )|$ $|\omega ^{n+1}|$ and (\ref{H7}) is reduced to%
\begin{equation}
dv_{\Sigma }=\pm \text{ }|\omega ^{n+1}|\text{ }\omega ^{1}\wedge \omega
^{2}\wedge ...\wedge \omega ^{n}.  \label{H11}
\end{equation}

\noindent For an (oriented) Riemannian manifold $M,$ we take $dv_{M}$ to be
the associated volume form. Then we can take $\omega ^{j}$'s to be an
orthonormal basis in (\ref{H11}). Hence $|\omega ^{n+1}|$ $=$ $1$ and $%
dv_{\Sigma }$ (up to sign) is nothing but the area form with respect to the
induced metric.

\bigskip

\textbf{Example A.4.} For $M$ being a pseudohermitian $3$-manifold, let $%
e_{1}$ $\in $ $T\Sigma $ $\cap $ $\xi $ denote the characteristic field on
the nonsingular domain in \cite{chmy}. Let $e_{2}$ $\equiv $ $Je_{1}$ and $%
\alpha $ denote a function such that $T$ $+$ $\alpha e_{2}$ $\in $ $T\Sigma
. $ Let $e^{1},$ $e^{2}$ (and $\Theta )$ be the coframe dual to $e_{1},$ $%
e_{2} $ (and $T).$ We can take $v_{1}$ $=$ $e_{1},$ $v_{2}$ $=$ ($T$ $+$ $%
\alpha e_{2})/\sqrt{1+\alpha ^{2}},$ and $v_{3}$ $=$ ($\alpha T-e_{2})/\sqrt{%
1+\alpha ^{2}}$ while $\omega ^{1}$ $=$ $e^{1},$ $\omega ^{2}$ $=$ ($\Theta
+\alpha e^{2})/\sqrt{1+\alpha ^{2}}$ and $\omega ^{3}$ $=$ ($\alpha \Theta
-e^{2})/\sqrt{1+\alpha ^{2}}.$ Note that $\omega ^{1}$ $\wedge $ $\omega
^{2} $ $\wedge $ $\omega ^{3}$ $=$ $e^{1}$ $\wedge $ $e^{2}$ $\wedge $ $%
\Theta $ is the standard volume form with respect to the adapted metric $h$ $%
\equiv $ $\Theta \otimes \Theta $ $+$ $\frac{1}{2}d\Theta (\cdot ,J\cdot )$.
Observe that $e^{1},$ $e^{2}$ are orthonormal with respect to the
semipositive inner product (\ref{H2}) since they are different from $\theta
^{1},$ $\theta ^{1^{\prime }}$ by an orthogonal transformation. Thus by (\ref%
{H2}) we have%
\begin{equation}
|\omega ^{3}|^{2}=\frac{<-e^{2},-e^{2}>}{1+\alpha ^{2}}=\frac{1}{1+\alpha
^{2}}  \label{H12}
\end{equation}

\noindent and 
\begin{eqnarray}
\omega ^{1}\wedge \omega ^{2} &=&e^{1}\wedge (\Theta +\alpha e^{2})/\sqrt{%
1+\alpha ^{2}}  \label{H13} \\
&=&\sqrt{1+\alpha ^{2}}e^{1}\wedge \Theta  \notag
\end{eqnarray}

\noindent on $\Sigma $ by noting that $e^{1}$ $\wedge $ $e^{2}$ $=$ $\alpha
e^{1}$ $\wedge $ $\Theta $ on $\Sigma .$ Substituting (\ref{H12}), (\ref{H13}%
) into (\ref{H11}) with $n=2,$ we conclude%
\begin{equation}
dv_{\Sigma }=\pm |\omega ^{3}|\omega ^{1}\wedge \omega ^{2}=\pm e^{1}\wedge
\Theta .  \label{H14}
\end{equation}

\noindent The above expression first appeared in \cite{chmy}.

\bigskip

\textbf{Example A.5. }Let $\pi _{\xi }^{h}:$ $TM\rightarrow $ $\xi $ denote
the projection onto $\xi $ according to the adapted metric $h.$ Then we have
($|\cdot |_{h}$ denotes the length with respect to $h)$%
\begin{equation}
|\pi _{\xi }^{h}(v_{3})|_{h}=|\frac{-e_{2}}{\sqrt{1+\alpha ^{2}}}|_{h}=\frac{%
1}{\sqrt{1+\alpha ^{2}}}=|\omega ^{3}|.  \label{H15}
\end{equation}

\noindent In view of (\ref{H15}), (\ref{H14}) and $v_{3}$ being a unit
normal with respect to $h$, we obtain%
\begin{equation}
dv_{\Sigma }=\pm |\pi _{\xi }^{h}(N)|_{h}d\Sigma _{h}  \label{H16}
\end{equation}

\noindent where $N$ denotes the unit normal (unique up to sign) with respect
to $h$ and $d\Sigma _{h}$ denotes the area element with respect to the
metric induced from $h.$ The expression (\ref{H16}) appeared in \cite{Ri2}
for $M$ being the $3$-dimensional Heisenberg group.

\bigskip

Next we are going to deduce a formula for the mean curvature $H$ viewed as
the first variation of the area. Recall that in (\ref{H7}) and (\ref{H6}) $%
\varphi $ is a defining function of a hypersurface $\Sigma $ in a manifold $%
M $ of dimension $n+1$ and $\omega ^{j}$'s are independent 1-forms. We
assume further $\omega ^{n+1}$ $=$ $0$ on $\Sigma $ and 
\begin{equation}
d\omega ^{j}=\sum_{k=1}^{n+1}\omega ^{k}\wedge \omega _{k}^{j},\text{ }%
j=1,...,n+1  \label{H17}
\end{equation}

\noindent for some 1-forms $\omega _{k}^{j}.$ Starting from (\ref{H7}), we
compute%
\begin{eqnarray}
&&\delta _{fv_{n+1}}\int_{\Sigma }dv_{\Sigma }  \label{H18} \\
&=&\int_{\Sigma }L_{fv_{n+1}}\{\frac{(-1)^{n}}{v_{n+1}(\varphi )}|d\varphi
|\omega ^{1}\wedge \omega ^{2}\wedge ...\wedge \omega ^{n}\}  \notag
\end{eqnarray}

\noindent where $f$ is a $C^{\infty }$ smooth function and $L_{fv_{n+1}}$
denotes the Lie derivative in the direction $fv_{n+1}.$ Let $i_{X}\eta $
denote the interior product of the vector field $X$ and the differential
form $\eta .$ Observe that $i_{fv_{n+1}}(\omega ^{1}\wedge \omega ^{2}\wedge
...\wedge \omega ^{n})$ $=$ $0$ since $\omega ^{j}(v_{n+1})$ $=$ $0$ for $1$ 
$\leq $ $j$ $\leq $ $n.$ It follows from $L_{fv_{n+1}}$ $=$ $%
i_{fv_{n+1}}\circ d$ $+$ $d\circ i_{fv_{n+1}}$ that%
\begin{eqnarray}
&&L_{fv_{n+1}}\{\frac{(-1)^{n}}{v_{n+1}(\varphi )}|d\varphi |\omega
^{1}\wedge \omega ^{2}\wedge ...\wedge \omega ^{n}\}  \label{H19} \\
&=&i_{fv_{n+1}}\circ d\{\frac{(-1)^{n}}{v_{n+1}(\varphi )}|d\varphi |\omega
^{1}\wedge \omega ^{2}\wedge ...\wedge \omega ^{n}\}.  \notag
\end{eqnarray}

\noindent Compute%
\begin{eqnarray}
&&d(\omega ^{1}\wedge \omega ^{2}\wedge ...\wedge \omega ^{n})  \label{H20}
\\
&=&\sum_{j=1}^{n}(-1)^{j-1}\omega ^{1}\wedge ...\wedge d\omega ^{j}\wedge
...\wedge \omega ^{n}  \notag \\
&=&-(\sum_{j=1}^{n}\omega _{j}^{j})\wedge \omega ^{1}\wedge \omega
^{2}\wedge ...\wedge \omega ^{n}  \notag \\
&&+(-1)^{n}\sum_{j=1}^{n}\omega ^{1}\wedge ...\wedge \omega _{n+1}^{j}\wedge
...\wedge \omega ^{n}\wedge \omega ^{n+1}  \notag
\end{eqnarray}

\noindent where we have used (\ref{H17}). We can now obtain from (\ref{H18}%
), (\ref{H19}), and (\ref{H20}) that%
\begin{eqnarray}
&&\delta _{fv_{n+1}}\int_{\Sigma }dv_{\Sigma }  \label{H21} \\
&=&\int_{\Sigma }f\{v_{n+1}(\frac{|d\varphi |}{v_{n+1}(\varphi )})  \notag \\
&&+(-1)^{n}\frac{|d\varphi |}{v_{n+1}(\varphi )}\sum_{j=1}^{n}(\omega
_{n+1}^{j}(v_{j})-\omega _{j}^{j}(v_{n+1}))\}\omega ^{1}\wedge \omega
^{2}\wedge ...\wedge \omega ^{n}.  \notag
\end{eqnarray}

\noindent Here we have used $\omega ^{n+1}$ $=$ $0$ on $\Sigma .$ Comparing (%
\ref{H21}) with (\ref{H7}) we obtain the mean curvature%
\begin{eqnarray}
H &=&\mp \{(-1)^{n}\frac{v_{n+1}(\varphi )}{|d\varphi |}v_{n+1}(\frac{%
|d\varphi |}{v_{n+1}(\varphi )})  \label{H22} \\
&&+\sum_{j=1}^{n}(\omega _{n+1}^{j}(v_{j})-\omega _{j}^{j}(v_{n+1}))\}. 
\notag
\end{eqnarray}

\bigskip

\textbf{Example A.6.} In the Riemannian case, we can take $\omega ^{j},$ $%
\omega _{k}^{j}$ in (\ref{H17}) to be an orthonormal coframe and the
associated connection forms, resp. such that $\omega ^{n+1}$ $=$ $0$ on $%
\Sigma .$ So from $<\omega ^{i},\omega ^{j}>$ $=$ $\delta _{ij}$ and $%
v_{i}(\varphi )$ $=$ $0$ for $1$ $\leq $ $i$ $\leq $ $n$ we have%
\begin{eqnarray*}
|d\varphi |^{2} &=&\sum_{i,j=1}^{n+1}v_{i}(\varphi )<\omega ^{i},\omega
^{j}>v_{j}(\varphi ) \\
&=&\sum_{i=1}^{n+1}(v_{i}(\varphi ))^{2}=(v_{n+1}(\varphi ))^{2}.
\end{eqnarray*}%
\noindent It follows that 
\begin{equation}
\frac{|d\varphi |}{v_{n+1}(\varphi )}=\pm 1.  \label{H23}
\end{equation}%
\noindent On the other hand, the Riemannian connection forms $\omega
_{k}^{j} $'s satisfy the skew-symmetric condition: $\omega _{k}^{j}$ $+$ $%
\omega _{j}^{k}$ $=$ $0.$ So we have%
\begin{equation}
\omega _{j}^{j}=0  \label{H23.1}
\end{equation}%
\noindent Write $\omega _{i}^{n+1}$ $=$ $h_{ij}\omega ^{j}$ where $h_{ij}$ $%
= $ $h_{ji}$ (due to $d\omega ^{n+1}$ $=$ $0$ on $\Sigma )$ are known to be
coefficients of the second fundamental form. We then have%
\begin{eqnarray}
\sum_{j=1}^{n}\omega _{n+1}^{j}(v_{j}) &=&-\sum_{j=1}^{n}\omega
_{j}^{n+1}(v_{j})  \label{H24} \\
&=&-\sum_{j=1}^{n}h_{jj}  \notag
\end{eqnarray}

\noindent Substituting (\ref{H23}), (\ref{H23.1}), and (\ref{H24}) into (\ref%
{H22}), we obtain%
\begin{equation*}
H=\pm \sum_{j=1}^{n}h_{jj}.
\end{equation*}

\noindent This verifies the formula (\ref{H22}) for the Riemannian situation.

\bigskip

\textbf{Example A.7.} Consider a surface $\Sigma $ in a pseudohermitian
3-manifold. We will continue to use the notations in Example 3.3. Take $%
\omega ^{1}$ $=$ $\Theta ,$ $\omega ^{2}$ $=$ $e^{1},$ and $\omega ^{3}$ $=$ 
$e^{2}$ $-$ $\alpha \Theta .$ That $\omega ^{3}$ $=$ $0$ on $\Sigma $
follows from $e_{1}$ $\in $ $T\Sigma $ and $T$ $+$ $\alpha e_{2}$ $\in $ $%
T\Sigma .$ The corresponding dual vectors are $v_{1}$ $=$ $T$ $+$ $\alpha
e_{2},$ $v_{2}$ $=$ $e_{1},$ and $v_{3}$ $=$ $e_{2}.$ Since $<\omega
^{i},\omega ^{j}>$ $=$ $\delta _{ij}$ by (\ref{H2}) and $v_{i}(\varphi )$ $=$
$0$ for $i$ $=$ $1,$ $2,$ we still have (\ref{H23}) with $n+1$ $=$ $3.$ Here 
$\varphi $ is a defining function of $\Sigma .$ From the structure equations
(A.1r), (A.3r) in \cite{chmy}, we can take%
\begin{eqnarray}
\omega _{1}^{1} &=&0,\text{ }\omega _{2}^{1}=2e^{2},\text{ }\omega _{3}^{1}=0
\label{H25} \\
\omega _{2}^{2} &=&0,\text{ }\omega _{3}^{2}=-\omega  \notag
\end{eqnarray}

\noindent in (\ref{H17}), where $i\omega $ is the pseudohermitian connection
form. Now by (\ref{H22}), (\ref{H23}) and (\ref{H25}), we have%
\begin{eqnarray*}
H &=&\mp (0+\sum_{j=1}^{2}(\omega _{3}^{j}(v_{j})-\omega _{j}^{j}(v_{3})) \\
&=&\mp \omega _{3}^{2}(e_{1})=\pm \omega (e_{1}).
\end{eqnarray*}

\noindent This is an expression of the $p$-mean curvature in \cite{chmy}.

\end{document}